\input amstex
\documentstyle{amsppt}
\NoRunningHeads
\NoBlackBoxes
\document

\def\Q{{\Bbb Q}}
\def\JJ{{\Bbb J}}

\def\ell{{\text{ell}}}

\def\RR{{\Bbb R}}

\def\h1{\hat{\bold 1}}

\def\Hom{\text{Hom}}

\def\Ua{U_q(\tilde\g)}
\def\U2{{\Ua}_2}
\def\g{\frak g}
\def\n{\frak n}

\def\Z{\Bbb Z}
\def\C{\Bbb C}
\def\d{\partial}

\def\Tr{\text{\rm Tr}}
\def\l{\lambda}

\def\<{\langle}
\def\>{\rangle}
\def\o{\otimes}

\def\End{\text{End}}

\def\h{{\frak h}}
\topmatter
\title Traces of intertwiners for quantum groups
and difference equations, I
\endtitle

\author {\rm {\bf Pavel Etingof and Alexander Varchenko} \linebreak
\vskip .05in
Department of Mathematics, Rm 2-165\linebreak
MIT, 77 Mass. Ave.\linebreak 
Cambridge, MA 02139, USA\linebreak
e-mail: etingof\@math.harvard.edu
\vskip .03in
and
\vskip .03in
Department of Mathematics\linebreak
University of North Carolina\linebreak 
Chapel Hill, NC 27599, USA\linebreak
e-mail: av\@math.unc.edu}
\endauthor
\endtopmatter

\head 0. Introduction\endhead

This paper begins a series of papers whose 
goal is to establish a representation-theoretic interpretation 
of the quantum Knizhnik-Zamolodchikov-Bernard (qKZB) 
equations, and use this interpretation to study solutions 
of these equations. It was motivated by the recent work
on the qKZB equations \cite{F,FTV1-2,MV,FV2-5}, and by the theory 
of ``quantum conformal blocks'' that began with the classical paper \cite{FR}. 

{\bf 0.1.}  
The qKZB equations
\cite{F} are difference equations with respect to an unknown function
$f(z_1,...,z_N,\l,\tau,\mu,p)$ with values in 
$V_1\o...\o V_N\o V_N^*\o...\o V_1^*$, where $V_i$ are 
suitable finite dimensional representations of
the quantum group $U_q(\g)$ ($\g$ is a simple Lie algebra),  
$z_i,p,\tau\in \C$, and $\l,\mu$ are weights for $\g$. 

The qKZB equations 
are a q-deformation of the Knizhnik-Zamolodchikov-Bernard
(KZB) differential 
equations, and an
elliptic analogue of the quantum Knizhnik-Zamolodchikov 
(qKZ) difference equations, which are, in turn, generalizations 
of the usual (trigonometric) Knizhnik-Zamolodchikov equations. 

It is proved in \cite{FTV2} (using an integral representation of solutions) 
that for $\g=sl_2$ the monodromy of the qKZB equations
is given by the dual qKZB equations, which are obtained 
from the qKZB equations by interchanging $(\l,\tau)$ with $(\mu,p)$.
 This fact generalizes the
 monodromy theorems for the KZB and qKZ equations:
the monodromy of the KZB differential equations is the 
trigonometric degeneration 
of the qKZB equations (which involves 
dynamical R-matrices without spectral parameter) (see e.g. \cite{K}), 
and the monodromy of the qKZ equations
is given by elliptic dynamical R-matrices 
(but there is no difference equation) (\cite{TV1,TV2}; see also \cite{FR}).   

The self-duality of the qKZB equations leads one to expect that 
they should have symmetric solutions $u_{V_1,...,V_N}(z,\l,\tau,\mu,p)$, 
i.e. such that $u_{V_1,...,V_N}(z,\l,\tau,\mu,p)=
u^*_{V_N^*,...,V_1^*}(z,\mu,p,\lambda,\tau)$, where 
$u^*$ is the dual of $u$ (considered as an endomorphism 
of $V_1\o...\o V_N$). 
Such a solution $u$ (for $\g=sl_2$) 
was constructed in \cite{FV2,FTV2}, by an explicit integral formula.
It is called the universal hypergeometric function.  
This function has many 
interesting properties, in particular the $SL(3,\Z)$-symmetry \cite{FV3-5}, 
where the group $SL(3,\Z)$ acts on the lattice 
$\Z^3$ generated by the periods $1,\tau,p$. 
A consequence of this symmetry is the qKZB-heat equation \cite{FV2}
for the function $u$, 
which is a q-deformation of the KZB-heat equation \cite{Ber}.  

{\bf 0.2.}
A central fact about the KZB and qKZ equations 
(and one of the main reasons why they are interesting) is that
they are satisfied by conformal blocks. 
More precisely, the KZB equations
are satisfied by conformal blocks of the Wess-Zumino-Witten 
conformal field theory on an elliptic curve \cite{Ber}, and the qKZ  
equations are satisfied by quantum conformal blocks on the
cylinder \cite{FR}. In representation theoretic terms, 
conformal blocks on an elliptic curve are traces 
of products of intertwining operators for affine Lie algebras
(weighted by an element from the maximal torus) 
\cite{Ber}, and quantum conformal blocks on the 
cylinder are highest matrix elements of products of intertwining operators  
for quantum affine algebras\cite{FR}. 
This representation theoretic interpretation of the KZB and qKZ equations  
is not only interesting by itself, but it also allows to 
prove nontrivial properties of solutions, e.g. monodromy theorems
(see e.g. \cite{K,FR}). 

The goal of this series
is to give a similar
interpretation of the qKZB equations. In light of the above, 
the main idea is obvious: one should consider 
quantum conformal blocks on an elliptic curve, or, 
representation theoretically, traces of products 
of intertwining operators for quantum affine 
algebras, weighted by an element of the maximal torus. It is natural to expect 
that such traces satisfy a pair of dual qKZB equations. 
This is actually true, and we plan to give a proof of it 
in a subsequent part of the series. However, the details of the proof are 
relatively complicated, and we would like to start with 
a simpler (``trigonometric'') limiting case, when $\tau,p\to\infty$. This
limiting case is the main subject of this paper. 

{\bf 0.3.} The structure of this paper is as follows. 

In Section 1 we introduce the main object of the paper -- the 
renormalized universal trace function
$F_{V_1,...,V_N}(\l,\mu)\in (V_1\o...\o V_N)[0]\o (V_N^*\o...\o V_1^*)[0]$, 
where $\l,\mu$ are weights for $\g$. 
It is obtained from traces of products of intertwining 
operators for $U_q(\g)$ weighted by an element of the maximal torus. 
At the end of the section we formulate the main results of the 
paper -- Theorems 1.1--1.5. 

Theorems 1.1 and 1.2 state that the function $F_{V_1,...,V_N}$ 
satisfies two systems of difference equations, one with shifts of $\l$, 
and the other with shifts of $\mu$, which go to each other under 
the transformation $\l\to \mu$, $\mu\to \l$. In the special case $\g=sl_n$, 
$N=1$, $V_1=S^{mn}\C^n$, these systems (as was shown in \cite{EK1})
reduce to the 
trigonometric Macdonald-Ruijsenaars (MR) systems, so we call them 
the MR system and the dual MR system.   

Theorem 1.5 states that the function $F_{V_1,...,V_N}(\l,\mu)$ is symmetric: 
$F_{V_1,...,V_N}(\l,\mu)=F_{V_N^*,...,V_1^*}^*(\mu,\l)$,
where $*$ is the permutation of components.
It follows from Theorems 1.1 and 1.2. 

Theorems 1.3 and 1.4 state that the function $F_{V_1,...,V_n}$ 
satisfies two additional systems of difference equations --
the trigonometric degenerations of the 
qKZB and the dual qKZB equations, respectively. 

Theorems 1.1-1.5 are proved in Sections 2-5.

In Section 6, we study the symmetry of trace functions under $q\to q^{-1}$, 
and define a modified trace function $u_V(\l,\mu)$, by renormalizeing 
$F_V(\l,\mu)$. (This function is introduced to connect our paper 
with the papers \cite{FV2-FV5}; here we define it only for $N=1$, and plan 
to define it in general in another paper).  
Using the $q\to q^{-1}$ transformation properties and Theorem 1.5, we 
show that the function $u_V$ is symmetric. 

In Section 7, we compute the function 
$F_V(\l,\mu)$, $u_V(\l,\mu)$ explicitly in the case $\g=sl_2$.

In Section 8, we compute explicitly the trigonometric degeneration of the 
function $u$ from \cite{FV2}, in the case $N=1$.
 We show that this function is the same as $u_V(\l,\mu)$ up to 
normalization.   

In Section 9, we explain that Macdonald's theory for root systems of type 
$A_{n-1}$  is a special case of the theory developed in this paper, 
for  $\g=sl_n$, 
$N=1$, $V_1=S^{mn}\C^n$.

In Section 10, we consider limiting (degenerate) cases of the theory 
developed in this paper. 

{\bf 0.4.} In subsequent papers of the series, we plan: 

1. To give a representation theoretic-proof of the qKZB heat equation
and the orthogonality relations for the 
trigonometric degeneration of the function $u$ (\cite{FV2}), 
using the ideas of \cite{EK1,EK2,EK3}. Cherednik's theory 
of difference Fourier transform and Macdonald-Mehta identities 
for root systems of type A is a special 
case of this theory, corresponding to the situation $\g=sl_n$, 
$N=1$, $V_1=S^{mn}\C^n$. 

2. To give a representation theoretic derivation of the resonance relations 
from \cite{FV3} in the trigonometric case, using the ideas of 
\cite{ES}. 

3. To generalize all the results to the case of quantum affine algebras.
This involves a representation theoretic definition 
of the function $u$ from \cite{FV2} for generic values of parameters, 
for any simple Lie algebra and representations, and a representation 
theoretic proof of its main properties, such as 
qKZB and MR equations, orthogonality, modular transformations
(e.g. the qKZB heat equation).
As a special case, this theory should contain Macdonald's theory 
for affine root systems of type $\hat A_{n-1}$, which was originated in 
\cite{EK4} but has not been developed from an analytic standpoint. 
In particular, the classical limit ($q\to 1$) of the modular 
transformation of the function $u$ should yield the result 
of Kirillov \cite{K1,K2} which says that the modular transformation $S$ 
of affine Jack polynomials (which are essentially the 1-point functions 
of the WZW model in genus 1, see \cite{EK4}) is given 
by a matrix of special values Macdonald polynomials at roots of unity.    

4. Specializing this theory to the critical level, to 
prove that radial parts of the central elements 
of $U_q(\widehat{sl_n})$ at the critical level corresponding to 
the representation $S^{mn}\C^n$, are elliptic Ruijsenaars 
operators (as far as we know, this is known only in the trigonometric 
degeneration). 

{\bf Acknowledgments.} The first author was partially supported by 
the NSF grant DMS-9700477, and
thanks the UNC mathematics department for hospitality.
The work of the first author was partly done while he was employed
by the Clay Mathematics Institute as a CMI Prize Fellow.  
The second author was supported 
by the NSF grant DMS-9801582, and 
is grateful to the Harvard 
Mathematics department for hospitality. 
The authors are grateful to A.Kirillov Jr. for useful 
suggestions on how to improve the paper. 

\head 1. Trace functions for $U_q(\g)$\endhead

\subhead 1.1. The trace functions\endsubhead

Let $\g$ be a simple Lie algebra over $\C$. Let 
$\h$ be a Cartan subalgebra of $\g$, and $\alpha_i$ be simple roots of $\g$, 
$i=1,...,r$. Let $(a_{ij})$ be the Cartan matrix of $\g$. 
Let $d_i$ be relatively prime positive integers such that $(d_ia_{ij})$ 
is a symmetric matrix. Let 
$e_i,f_i,h_i$ be the Chevalley generators of $\g$. 

Let $t$ be a complex number which is not purely imaginary, 
and $q=e^t$. For any operator $A$, we will denote $e^{tA}$ by $q^A$. 

Let $U_q(\g)$  be the Drinfeld-Jimbo quantum group corresponding 
to $\g$. We will use the same definition 
of $U_q(\g)$ as in \cite{EFK}. Namely $U_q(\g)$ is a Hopf algebra 
with generators
$E_i,F_i$, $i=1,\dots,r$, $q^h$, $h\in\h$, with relations:
$$q^{x+y}=q^xq^y,\ x,y\in \h\quad 
q^hE_jq^{-h}=q^{\alpha_j(h)}E_i,\quad q^hF_jq^{-h}=q^{-\alpha_j(h)}F_i$$
$$E_iF_j-F_jE_i=\delta_{ij}\frac{q^{d_ih_i}-q^{-d_ih_i}}{q^{d_i}-q^{-d_i}},$$
$$\sum_{k=0}^{1-a_{ij}} (-1)^k \bmatrix 1-a_{ij} \\ k \endbmatrix_{q_i}
E_i^{1-a_{ij}-k}E_jE_i^k=0,\quad i\ne j,$$
$$\sum_{k=0}^{1-a_{ij}} (-1)^k \bmatrix 1-a_{ij} \\ k \endbmatrix_{q_i}
F_i^{1-a_{ij}-k}F_jF_i^k=0,\quad i\ne j.$$
where $q_i = q^{d_i}$ and we used notation
$$ \bmatrix n \\ k \endbmatrix _q=\frac{[n]_q!}{[k]_q! [n-k]_q!},
\quad  [n]_q! = [1]_q \cdot [2]_q \cdot \dots \cdot [n]_q,
\quad  [n]_q= \frac {q^n - q^{-n}}{q-q^{-1}}$$

Comultiplication $\Delta,$ antipode $S$ and counit $\epsilon$
in $U_q(\g)$ are given by
$$\Delta(E_i) = E_i\otimes q^{d_ih_i} + 1\otimes E_i, \quad
  \Delta(F_i) = F_i\otimes 1 + q^{-d_ih_i}\otimes F_i, \quad
  \Delta(q^h) = q^h \otimes q^h$$
$$S(E_i)=-E_iq^{-d_ih_i},\quad S(F_i)=-q^{d_ih_i}F_i,\quad S(q^h)=q^{-h}$$
$$\epsilon(E_i) = \epsilon(F_i) = 0,\quad \epsilon(q^h) = 1$$

Let $M_\mu$ be the Verma module over $U_q(\g)$ with highest weight 
$\mu$, $v_\mu$ its highest weight vector. 
Let $V$ be a finite dimensional representation of $U_q(\g)$, 
and $v\in V$ a vector of weight $\mu_v$. It is well known 
that for a generic $\mu$ there exists a unique intertwining 
operator $\Phi_\mu^v:M_\mu\to M_{\mu-\mu_v}\o V$ such that 
$\Phi_\mu^vv_\mu=v_{\mu-\mu_v}\o v+l.o.t.$ (here l.o.t. denotes lower order 
terms). It is useful to consider the generating function 
of such operators, $\Phi_\mu^V\in \Hom_\C(M_\mu,\oplus M_\nu\o V\o V^*)$, 
defined by $\Phi_\mu^V=\sum_{v\in B}\Phi_\mu^v\o v^*$, where the summation is 
taken over a homogeneous basis $B$ of $V$. 

Let $V_1,...,V_N$ be finite dimensional representations of $U_q(\g)$, 
$v_i\in V_i$ vectors of weights $\mu_{v_i}$, 
such that $\sum\mu_{v_i}=0$. 
Define the following formal series
in \linebreak $V_1\o...\o V_N[0]\o q^{2(\l,\mu)}\C 
[[q^{-2(\l,\alpha_1)},...,q^{-2(\l,\alpha_r)}]]$: 
$$
\Psi^{v_1,...,v_N}(\l,\mu)=
\text{Tr}|_{M_\mu}
((\Phi^{v_1}_{\mu-\sum_{i=2}^N\mu_{v_i}}\o 1^{N-1})
...\Phi^{v_N}_\mu q^{2\l}) 
\tag 1.1
$$
It follows from Corollary 3.4, \cite{ES} 
that this series converges (in a suitable region of values 
of the parameters) to a function
of the form $q^{2(\l,\mu)}f(\l,\mu)$, where $f$ is a rational function 
in $q^{2(\l,\alpha_i)}$ and $q^{2(\mu,\alpha_i)}$,
which is a finite sum of products of functions of $\l$ and
functions of $\mu$. The function (1.1) 
will be 
called {\it a trace function}. 

Define also the universal trace function, 
with values in $V_1\o...\o V_N\o V_N^*\o...\o V_1^*$:
$$
\Psi_{V_1...V_N}(\l,\mu)=\sum_{v_i\in B_i} 
\Psi^{v_1,...,v_N}(\l,\mu)\o v_N^*\o ...\o v_1^*,\tag 1.2
$$
where $B_i$ are homogeneous bases of $V_i$. 
It is easy to see that this function takes values in 
$(V_1\o..\o V_N)[0]\o (V_N^*\o...\o V_1^*)[0]$. 

Using the generating functions $\Phi_\mu^V$,
one can express the universal trace function as follows:
$$
\Psi_{V_1...V_N}(\l,\mu)=
\text{Tr}((\Phi^{V_1}_{\mu+\sum_{i=2}^Nh^{(*i)}}\o 1^{N-1})
...\Phi^{V_N}_\mu q^{2\l}),
\tag 1.3
$$
where we label the components $V_i$ by $i$ and $V_i^*$ by $*i$, and 
the notation $h^{(k)}$ for a label $k$ was defined in \cite{F}: when 
acting on a homogeneous multivector, 
$h^{(k)}$ has to be replaced with the weight 
in the k-th component. 

{\bf Example 1.} (See Section 7) 
Let $\g=sl_2$. In this case let us represent weights by complex
numbers, so that the unique fundamental weight corresponds to
$1$. 
If $N=1$, and $V=V_1$ is the 3-dimensional 
representation, then 
$$
\Psi_V(\l,\mu)=\frac{q^{\l\mu}}{1-q^{-2\l}}\left(
1+(q^2-q^{-2})\frac{q^{-2\l}}{(1-q^{2\mu})(1-q^{-2(\l-1)})}\right)
$$
(since $V[0]$ is 1-dimensional, we view 
the function $\Psi_V$ as a scalar function).
This example is also computed in \cite{ES1}. 

\subhead 1.2. The main results\endsubhead

It turns out that the trace function $\Psi_{V_1...V_N}(\l,\mu)$ 
satisfies some remarkable difference equations. 
These equations are written in terms of so called dynamical 
R-matrices. Below we give a brief introduction to the theory of dynamical
R-matrices, referring the reader to the expository paper
\cite{ESch}
for a more detailed discussion of them.
 
Let $V,W$ be finite dimensional representations of $U_q(\g)$. 
Recall from \cite{EV} the definitions of the fusion matrix and the 
exchange matrix.

The fusion matrix is the operator $J_{WV}(\mu):W\o V\to W\o V$ defined 
by the formula
$$
(\Phi_{\mu-\mu_v}^{w}\o 1)\Phi_{\mu}^{v}=
\Phi_\mu^{J_{WV}(\mu)(w\o v)}.\tag 1.4
$$
The exchange matrix
$R_{VW}(\mu)\in \End(V\o W)$ is defined by 
$$
R_{VW}(\mu)=J_{VW}^{-1}(\mu)\Cal R^{21}|_{V\o W}J_{WV}^{21}(\mu),
\tag 1.5
$$ 
where 
$\Cal R$ is the universal R-matrix of $U_q(\g)$. 

We will also use the universal 
fusion matrix $J(\l)$ and the universal exchange matrix $R(\l)$. 
They take values in a completion of 
$U_q(\g)\o U_q(\g)$ and are defined by the condition that 
they give $J_{VW}(\l)$, 
$R_{VW}(\l)$ when evaluated in representations $V,W$
(cf. also \cite{ABRR},\cite{JKOS}). 

{\bf Remark.} 
The fusion matrix 
describes how to ``fuse'' together two intertwining operators. 
The exchange matrix describes how to exchange the order 
of intertwining operators. 
The fusion matrix satisfies the 2-cocycle identity (see below),
and is sometimes referred to as a ``quasi-Hopf twist''. The 
exchange matrix satisfies the quantum dynamical Yang-Baxter equation,  
and is refereed to as a ``dynamical R-matrix''. 

Let $\JJ (\l):=J(-\l-\rho)$, where $\rho$ is the half-sum of positive 
roots, and let $\RR (\l)=R(-\l-\rho)$. Let
$\Q(\l)=m_{21}(1\o S^{-1})(\JJ(\l))$,
where $m_{21}(a\o b):=ba$, and $S$ is the antipode. 
It is easy to show that $\Q(\l)$ is invertible for generic $\l$. 
  
Define 
$$
\JJ ^{1...N}(\l)=\JJ ^{1,2...N}(\l)\JJ ^{2,3...N}(\l)...
\JJ ^{N-1,N}(\l),\tag 1.6
$$
where, for example, $\JJ^{1,2...N}$ stands for 
$(1\o\Delta_{n-1})(\JJ)$, where $\Delta_{n-1}:U_q(\g)\to 
U_q(\g)^{n-1}$ is the iterated coproduct. We agree 
that $\JJ^1(\l)=1$. Thus $\JJ^{1...N}(\l)$ describes how to ``fuse'' 
$N$ intertwining operators.  

Let 
$$
\delta_q(\l)=q^{2(\l,\rho)}\prod_{\alpha>0}(1-q^{-2(\l,\alpha)})\tag 1.7
$$
 be the Weyl denominator.

 Set 
 $$
 \varphi_{V_1...V_N}(\l,\mu)=\JJ ^{1...N}(\l)^{-1}
 \Psi_{V_1...V_N}(\l,\mu) \delta_q(\l).\tag 1.8
 $$

 Finally, introduce the renormalized trace function
 $$
 F_{V_1...V_n}(\l,\mu)=
 [\Q^{-1}(\mu)^{(*N)}\o ...\o \Q^{-1}(\mu-h^{(*2...*N)})^{(*1)}]
 \varphi_{V_1...V_N}(\l,-\mu-\rho).\tag 1.9
 $$ 
 It is convenient to 
 formulate properties of trace functions using this renormalization.

 {\bf Example 2.} (See Section 7)
 If $\g=sl_2$, $N=1$, and $V=V_1$ is the 3-dimensional 
 representation, then 
 $$
 F_V(\l,\mu)=
 q^{-\l\mu}\frac{q^{2(\l+\mu)}-q^{2\l-2}-q^{2\mu-2}+1}{(1-q^{2\l-2})
 (1-q^{2\mu-2})}.
 $$
 (This formula is obtained after simplifications from formula 
 (7.20) below when $m=1$).  
 Note that it is seen from this formula that $F_V$ is 
 symmetric in $\l$ and $\mu$. 

 The following theorems describe the properties of $F_{V_1...V_N}(\l,\mu)$.

 For any finite dimensional $U_q(\g)$-module $W$, define 
 the difference operator $\Cal D_W$ acting on functions of 
 $\l\in \h^*$ 
 with values in $(V_1\o...\o V_N)[0]$ given by the formula
 $$
 \Cal D_W=\sum_\nu\Tr|_{W[\nu]}(\RR _{WV_1}^{01}(\l+h^{(2...N)})...
 \RR_{WV_N}^{0N}(\l))
 T_\nu, \tag 1.10
 $$
 where $T_\nu f(\l)=f(\l+\nu)$, and the component $W$ is labeled by $0$.

 \proclaim{Theorem 1.1} (Macdonald-Ruijsenaars equations)
 $$
 \Cal D_W^\l F_{V_1...V_N}(\l,\mu)=\chi_W(q^{-2\mu})
 F_{V_1...V_N}(\l,\mu),
 \tag 1.11
 $$
 where $\chi_W(x)=\sum \dim W[\nu]x^\nu$ is the character of $W$,
 and by $\Cal D_W^\l$ we mean the operator $\Cal D_W$ acting on $F$ as a function 
 of $\l$, in components $V_1,...,V_N$.
 \endproclaim

 Theorem 1.1 is proved in Section 2. 

 {\bf Example 3.}
 If $\g=sl_2$, $N=1$, and $V=V_1$ is the 3-dimensional 
 representation, and $W$ is the 2-dimensional representation, then 
 $$
 \Cal D_W=T+\frac{(1-q^{2\l-4})(1-q^{2\l+2})}
 {(1-q^{2\l-2})(1-q^{2\l})}T^{-1},
 $$
 where $Tf(\l)=f(\l+1)$.
 To prove this, it is enough to check that this operator is the unique operator 
 of the form $a(\l)T+b(\l)T^{-1}$ such that (1.11) holds for the function $F_V$ 
 given above (note that in our case $\chi_W(q^{-2\mu})=q^\mu+q^{-\mu}$). 

 Introduce also the dual Macdonald-Ruijsenaars operators 
 $\Cal D_W^\vee$, 
 acting on functions of 
 $\mu\in \h^*$ 
 with values in $(V_N^*\o...\o V_1^*)[0]$, by the formula
 $$
 \Cal D_W^\vee=\sum_\nu\Tr|_{W[\nu]}(\RR _{WV_N^*}^{01}(\mu+h^{(*1...*N-1)})...
 \RR_{WV_1^*}^{0N}(\mu))
 T_\nu, \tag 1.12
 $$
 where $V_j^*$ is considered as a module over $U_q(\g)$ via the antipode.

 \proclaim{Theorem 1.2} (Dual Macdonald-Ruijsenaars equations)
 $$
 \Cal D_W^{\vee,\mu} F_{V_1...V_N}(\l,\mu)=\chi_W(q^{-2\l})
 F_{V_1...V_N}(\l,\mu),
 \tag 1.13
 $$
 where $\Cal D_W^{\vee,\mu}$ is $\Cal D_W^\vee$ acting on $F$ as a function of $\mu$, 
 in components $V_N^*,...,V_1^*$.
 \endproclaim

 Theorem 1.2 is proved in Section 3. 

 To formulate the next two results, we need to define 
 some operators acting on functions of 
 $\l$ and $\mu$ with values in $(V_1\o...\o V_N)[0]\o 
 (V_N^*\o...\o V_1^*)[0]$. 

 For $j=1,...,N$ define the operators 
 $$
 D_j=(q^{-2\mu-\sum x_i^2})_{*j}(q^{-2\sum x_i\o x_i})_{*j,*1...*j-1},\tag 1.14
 $$
 where $x_i$ is an orthonormal basis of $\h$. Also, define the operators
 $$
 \gather
 K_j=\\
 \RR _{j+1,j}(\l+h^{(j+2...N)})^{-1}...\RR _{Nj}(\l)^{-1}
 \Gamma_j\RR _{j1}(\l+h^{(2...j-1)}+h^{(j+1...N)})...\RR _{jj-1}
 (\l+h^{(j+1...N)}),\tag 1.15\endgather
 $$
 where $\Gamma_jf(\l):=f(\l+h^{(j)})$, and $h^{j...k}$
 acting on a 
 homogeneous multivector has to be replaced with the sum of weights of 
 components $j,...,k$ of this multivector.
 It is easy to check that $D_j$ commute with each other, 
 and its is known that so do $K_j$ (\cite{F}).

{\bf Remark.} As we have mentioned before, these operators are 
the trigonometric limits of the qKZB operators with spectral
parameters, introduced by Felder.

 Analogously, define the operators 
 $$
 D_j^\vee=
 (q^{-2\l-\sum x_i^2})_{j}(q^{-2\sum x_i\o x_i})_{j,j+1...N},\tag 1.16
 $$
 and 
 $$
 \gather
 K_j^\vee=
 \RR _{*j-1,*j}(\mu+h^{(*1...*j-2)})^{-1}...\RR _{*1,*j}(\mu)^{-1}
 \Gamma_{*j}\times \\
 \RR _{*j,*N}(\mu+h^{(*j+1...*N-1)}+h^{(*1...*j-1)})...
 \RR _{*j,*j+1}
 (\mu+h^{(*1...*j-1)}),\tag 1.17\endgather
 $$
 where $\Gamma_{*j}f(\mu)=f(\mu+h^{*j})$. 
 Like $D_j,K_j$, the operators $D_j^\vee, K_j^\vee$ commute.

 \proclaim{Theorem 1.3} (the qKZB equations) 
 The function $F_{V_1...V_N}$ satisfies 
 the qKZB equations:
 $$
 F_{V_1...V_N}(\l,\mu)=
 (K_j\o D_j)F_{V_1...V_N}(\l,\mu).\tag 1.18
 $$
 \endproclaim 

 \proclaim{Theorem 1.4} (the dual qKZB equations)
 $$
 F_{V_1...V_N}(\l,\mu)=
 (D_j^\vee\o K_j^\vee)F_{V_1...V_N}(\l,\mu).\tag 1.19
 $$
 \endproclaim

 {\bf Example 4.} Let $\g=sl_2$, $N=2$, $V_1=V_2=\C^2$ with standard
 basis $v_+,v_-$.  
 In this case $V_1\otimes V_2[0]$ is 2-dimensional with basis 
 $v_+\otimes v_-,v_-\otimes v_+$, and the action of the 
 dynamical R-matrix in this basis is
 $$
 \RR(\lambda)=
\pmatrix 1&
q^{-1/2}\frac{q^{-1}-q}{q^{-2\lambda}-1}&\\ \frac{q^{-1}-q}{q^{2\lambda}-1}
&\frac{(q^{-2\lambda}-q^2)
(q^{-2\lambda}-q^{-2})}{(q^{-2\lambda}-1)^2}\endpmatrix.
$$
Therefore, if $F_{V_1,V_2}(\lambda,\mu)$ 
is represented by a 2 by 2 matrix with respect to the above
basis then the qKZB equation corresponding to $j=2$ has the form
$$
\gather
\pmatrix \frac{(q^{-2\lambda}-q^2)
(q^{-2\lambda}-q^{-2})}{(q^{-2\lambda}-1)^2}
&\frac{q^{-1}-q}{q^{-2\lambda}-1}&
\\ \frac{q^{-1}-q}{q^{2\lambda}-1}
&1\endpmatrix
\pmatrix F_{11}(\lambda,\mu) & F_{12}(\lambda,\mu)\\ 
F_{21}(\lambda,\mu) & F_{22}(\lambda,\mu)\endpmatrix
\pmatrix q^{\mu} &0\\ 0& q^{-\mu}\endpmatrix=\\
\pmatrix F_{11}(\lambda+1,\mu) & F_{12}(\lambda-1,\mu)\\ 
F_{21}(\lambda+1,\mu) & 
F_{22}(\lambda-1,\mu)\endpmatrix
\endgather
$$
Here for convenience we took the shift operator $\Gamma_2$
from the left side of the equation to the right side. 
 
{\bf Remark 1.} We should warn the reader 
that the term ``qKZB equations'' is normally 
used for equations which contain elliptic dynamical R-matrices with
spectral parameters,
and are difference equations with respect to these parameters
$z_1,...,z_N$ (see e.g.\cite{FTV1}). 
The equations we consider are a limiting case 
of the ``genuine'' qKZB equations, when the modular parameter 
$\tau$ goes to infinity, and the ratios of the spectral
parameters 
$z_j/z_{j+1}$ go to zero, 
with $e^{-2\pi \text{Im}\tau}<<|z_j/z_{j+1}|<<1$.
Namely, the equations considered here are the equations satisfied 
by the limit (if it exists) of a solution of the ``genuine'' qKZB equations 
in the described asymptotic zone. Throughout this paper, we will
abuse terminology and 
use the term ``qKZB equations'' to refer to this limiting case.  

{\bf Remark 2.} If $N=1$, equations (1.18) and (1.19) are trivial: 
(1.18) says that the $V$-component of $F_V$ has zero weight, and 
(1.19) says that the $V^*$-component of $F_V$ has zero weight. 

{\bf Remark 3.} It is not hard to show using the quantum
dynamical Yang-Baxter equation for the dynamical R-matrices 
that the Macdonald-Ruijsenaars operators
commute with the qKZB operators. Similarly, the dual
Macdonald-Ruijsenaars
operators commute with the dual qKZB operators.    

{\bf Remark 4.} Theorems 1.1, 1.3 have the following interpretation. 
Suppose $v_1,...,v_N$ are homogeneous vectors in 
$V_1,...,V_N$ 
of weights $\nu_1,...,\nu_N$, $\sum \nu_i=0$.  
Then the function $(F_{V_1,..,V_N}(\l,\mu),v_1\o...\o v_N)$ 
is an common eigenfunction of the operators $\Cal D_W^\lambda$
and $K_j$ 
with eigenvalues equal to $\chi_W(q^{-2\mu})$ and \linebreak $\Lambda_j(\mu)=
q^{-2(\mu,\nu_j)+(\nu_j,\nu_j)+2\sum_{i<j}(\nu_i,\nu_j)}$,
respectively.  
Thus the trace functions provide a solution of the problem of simultaneous
diagonalization of the Macdonald-Ruijsenaars operators 
$\Cal D_W^\lambda$ and qKZB operators $K_j$.  

{\bf Remark 5.} The problem to deduce equations of type (1.18) for trace 
functions (in the case of affine Lie algebras) was suggested to 
the first author in 1992 by his adviser 
Igor Frenkel, as a topic of his Ph.D. thesis. However, the first author
failed to solve this problem at that time,
partly because the adequate framework, 
the theory of dynamical R-matrices, 
 was not around yet.  
 
\proclaim{Theorem 1.5} (the symmetry identities) 
The function $F_{V_1...V_N}$ is symmetric: 
$$
F_{V_1...V_N}(\l,\mu)=F^*_{V_N^*...V_1^*}(\mu,\l),\tag 1.20
$$
where $F^*$ is the result of interchanging the factors
$(V_1\o...\o V_N)[0]$ and \linebreak
$(V_N^*\o...\o V_1^*)[0]$. 
\endproclaim

Theorem 1.4 is proved in Section 4. 
Theorem 1.5 is proved in Section 5, using Theorems 1.1 and 1.2. 
Theorem 1.5 and Theorem 1.4 obviously imply Theorem 1.3.

{\bf Remark 1.} Theorem 1.3 can also be derived directly, 
using the method of Frenkel-Reshetikhin of 
derivation of the Knizhnik-Zamolodchikov equations. However, 
this derivation is rather long, and we don't give it here.  

{\bf Remark 2.} In the special case $\g=sl_n$, $N=1$, $V_1=S^{mn}\C^n$, 
the function $F$ is closely related to the kernel of Cherednik's 
difference Fourier transform for $sl_n$ (\cite{C}, and the symmetry 
theorem above (Theorem 1.5) is closely related to Cherednik's
theorem on the symmetry of the difference Fourier transform \cite{C}.  

{\bf Remark 3.} In the special case of Remark 2, Theorem 1.5
was proved in \cite{ES1} (Theorem 5.6). 

{\bf Remark 4.} Theorems 1.1-1.5 can be generalized to the case when 
$\g$ is any symmetrizable Kac-Moody algebra, and $V_i$
highest weight modules over $U_q(\g)$. In this case, the
functions $\Psi,F$ make sense as formal power series, and 
(if $\g$ is not finite dimensional) the operators 
$\Cal D_W$, $\Cal D_W^\vee$ are infinite difference operators:
they are infinite sums of terms $f(\l)T_\nu$. However, one can show 
that these sums make sense as operators on power series. 
We plan to discuss this elsewhere.
 
\head 2. The Macdonald-Ruijsenaars (MR) equations \endhead

\subhead 2.1. Radial parts\endsubhead

Let $V$ be a finite dimensional $U_q(\g)$ module. Then we have the following 
proposition. 

\proclaim{Proposition 2.1} (i) For any  
element $X$ of $U_q(\g)$ there exists a unique 
difference operator $D_X$ acting on $V[0]$-valued functions, 
such that 
$$
\Tr(\Phi^V_\mu Xq^{2\l})=D_X\Tr(\Phi^V_\mu q^{2\l}).\tag 2.1
$$

(ii) If $X$ is central then $D_{XY}=D_YD_X$
for all $Y\in U_q(\g)$. In particular, if $X$, 
$Y$ are central then $D_XD_Y=D_YD_X$.  
\endproclaim

\demo{Proof} This is proved in Section 6 of \cite{EK1} (the assumption 
$g=gl_n$, which is made throughout \cite{EK1}, is not important 
for the proof of this result).  
$\square$\enddemo

We will call the operator $D_X$ the radial part of $X$.

Now recall the Drinfeld-Reshetikhin construction 
of central elements of $U_q(\g)$ (\cite{D,R}). 
In this construction, one defines elements 
$C_W$ corresponding to finite dimensional representations $W$ of $U_q(\g)$
by the formula 
$$
C_W=\text{Tr}|_W(1\o \pi_W)(\Cal R^{21}\Cal R(1\o q^{2\rho})).\tag 2.2
$$
The map $W\to C_W$ defines a homomorphism of the Grothendieck ring 
of the category of finite dimensional representations of $U_q(\g)$
to the center of $U_q(\g)$. 

Define difference operators $\Cal M_W:=D_{C_W}$.  
 
\proclaim{Proposition 2.2} (i) $\Cal M_W\Cal M_U=\Cal M_U\Cal M_W=\Cal 
M_{W\o U}$.

(ii) $\Cal M_W\Psi_V(\l,\mu)=\chi_W(q^{2(\mu+\rho)})\Psi_V(\l,\mu)$,
where $\chi_W(x)=\sum_\nu \text{dim}W[\nu]x^\nu$ is the character of $W$. 
\endproclaim

\demo{Proof} See \cite{EK1} (the proof for $gl(n)$ generalizes tautologically 
to other Lie algebras).
$\square$\enddemo

This proposition shows that if we can compute $\Cal M_W$ explicitly
then we will get a system of difference equations for $\Psi_V$. 

\subhead 2.2. The difference equations\endsubhead

Let $G(\l)=q^{-2\rho}\Q^{-1}(\l)S(\Q)(\l-h)$.

\proclaim{Theorem 2.3} For any $V[0]$-valued function 
on $\h^*$, 
$$
(\Cal M_Wf)(\l)=\sum_\nu\Tr|_{W[\nu]}(G(\l+h)\RR _{WV}(\l))
f(\l+\nu).\tag 2.3
$$
\endproclaim

The proof of Theorem 2.3 occupies Sections 2.3-2.6. 

\subhead 2.3. The defining property of $J$\endsubhead

Define  
$$
\Cal J(\l):=J(-\l-\rho+\frac{1}{2}(h^{(1)}+h^{(2)})).\tag 2.4
$$

\proclaim{Lemma 2.4} (The Arnaudon-Buffenoir-Ragoucy-Roche 
equation, \cite{ABRR}, see also \cite{JKOS})
Let $V,W$ be finite dimensional $U_q(\g)$-modules. 
Then 
$$
\Cal R^{21}(q^{2\l})_1
\Cal J(\l)=
\Cal J(\l)
q^{\sum x_i\o x_i}(q^{2\l})_1.\tag 2.5
$$
Moreover, this solution is unique among solutions 
of the form $1+N$, where $N$ has only summands with positive weight 
of the second component. 
\endproclaim

{\bf Remark.} As pointed out in \cite{EV}, this equation
is a limiting case of the quantum Knizhnik-Zamolodchikov equation 
of \cite{FR}, obtained when the affine quantum group degenerates 
into the finite dimensional quantum group. 

\demo{Proof} We recall (\cite{ABRR}, see also 
\cite{EV}, section 9) that $J_{WV}(\l)$ 
satisfies the equation 
$$
J_{WV}(\lambda)(q^{2(\lambda+\rho)-\sum x_i^2})_2=
\Cal R_{VW}^{21}q^{-\sum x_i\otimes  x_i}(q^{2(\lambda+\rho)-\sum x_i^2})_2
J_{WV}(\lambda),\tag 2.6
$$ 
Using the weight zero property of $J$, we get
$$
J_{WV}(\l)(q^{-2(\l+\rho)})_1(q^{-\sum x_i^2})_2=
\Cal R_{VW}^{21}(q^{-2(\l+\rho)})_1q^{-\sum (x_i\otimes 1+1\o x_i)(1\o x_i)}
J_{WV}(\l),\tag 2.7
$$
Now apply both sides of (2.7) to the subspace $W[\nu]\o V[\mu]$. 
Using the weight zero condition again, we get (on that subspace) 
$$
J_{WV}(\l)(q^{-2(\l+\rho)})_1 q^{\mu+\nu}_1(q^\nu)_2=
\Cal R_{VW}^{21}(q^{-2(\l+\rho)})_1(q^{\mu+\nu})_1
J_{WV}(\l).\tag 2.8
$$
Thus, replacing $\l$ with $\l+\frac{1}{2}(\mu+\nu)$, we get
$$
J_{WV}(\l+\frac{1}{2}(\mu+\nu))(q^{-2(\l+\rho)})_1 (q^\nu)_2=
\Cal R_{VW}^{21}(q^{-2(\l+\rho)})_1
J_{WV}(\l+\frac{1}{2}(\mu+\nu)),\tag 2.9
$$
which means that
$$
\Cal R^{21}_{VW}(q^{-2(\l+\rho)})_1
J_{WV}(\l+\frac{1}{2}(h^{(1)}+h^{(2)}))=
J_{WV}(\l+\frac{1}{2}(h^{(1)}+h^{(2)}))
q^{\sum x_i\o x_i}(q^{-2(\l+\rho)})_1,\tag 2.10
$$
This implies the first statement.
The second statement is straightforward, and also
follows from \cite{ABRR}. 
$\square$\enddemo
 
We will also need the cocycle identity for $\Cal J$.
To deduce it, recall 
the 2-cocycle identity for $J$ (see \cite{A,EV}):
$$
J^{12,3}(\l)J^{12}(\l-h^{(3)})=J^{1,23}(\l)J^{23}(\l).\tag 2.11
$$
Thus the cocycle identity for $\Cal J$ has the form
$$
\Cal J^{12,3}(\l)\Cal J^{12}(\l+\frac{1}{2}h^{(3)})=
\Cal J^{1,23}(\l)\Cal J^{23}(\l-\frac{1}{2}h^{(1)}).\tag 2.12
$$

\subhead 2.4. The function $Z_V$\endsubhead
Let $V$ be a finite-dimensional $U_q(\g)$-module. 
Consider the following function with values in $V\o V^*\o U_q(\g)$
(with components labeled as $1,*1,2$):
$$
Z_V(\l,\mu)=
\Tr_0(\Phi^{V,01}_\mu\Cal R^{20}q^{2\l}).\tag 2.13
$$

\proclaim{Lemma 2.5}
$$
Z_V(\l,\mu)=
\Cal R^{21}(q^{2\l})_1Z_V(\l,\mu).\tag 2.14
$$
\endproclaim

\demo{Proof}
Let us move the R-matrix $\Cal R^{20}$ around the trace. We get 
$$
\gather
Z_V(\l,\mu)=
\Tr_0(\Cal R^{21}\Cal R^{20}\Phi^{V,01}_\mu q^{2\l})=\\
\Cal R^{21}
\Tr_0(\Phi^{V,01}_\mu q^{2\l}\Cal R^{20}),\tag 2.15\endgather
$$
and the Lemma follows after interchangicng $q^{2\l}$ with $\Phi_\mu^V$ 
and moving it around the trace.  
$\square$\enddemo

\proclaim{Corollary 2.6} 
$$
Z_V(\l,\mu)=\Cal J^{12}(\l)
\Psi_V^{(1)}(\l+\frac{1}{2} h^{(2)},\mu).\tag 2.16
$$
\endproclaim 

\demo{Proof} Both sides of the equation satisfy (2.14) 
and have the form 
$\Psi_V^{(1)}(\l+\frac{1}{2} h^{(2)},\mu)+N$, where $N$ has positive 
degree in component 2. It is easy to show that a solution of (2.14) 
with such property is unique, which implies the lemma. 
$\square$\enddemo

\subhead 2.5. The function $X_V$\endsubhead 

Define the function 
$$
X_V(\l,\mu)=
\Tr_0(\Phi^{V,01}_\mu\Cal R^{20}q^{2\l} (\Cal R^{03})^{-1}),\tag 2.17
$$
with values in $V\o V^*\o U_q(\g)\o U_q(\g)$ (components labeled $1,*1,2,3$). 

\proclaim{Lemma 2.7} 
$$
X_V(\l,\mu)=
\Cal R^{12,3}(q^{2\l})_3
X_V(\l,\mu)
(q^{-2\l})_3(\Cal R^{23})^{-1}.\tag 2.18
$$
\endproclaim

\demo{Proof} 
We have  
$$
\gather
X_V(\l,\mu)=
\Tr_0((\Cal R^{03})^{-1}\Phi^{V,01}_\mu\Cal R^{20}q^{2\l})=\\
\Cal R^{13}\Tr_0(\Phi^{V,01}_\mu(\Cal R^{03})^{-1}
\Cal R^{20}q^{2\l})=\\
\Cal R^{13}\Cal R^{23}
\Tr_0(\Phi^{V,01}_\mu\Cal R^{20}(\Cal R^{03})^{-1}q^{2\l})
(\Cal R^{23})^{-1}=\\
\Cal R^{13}\Cal R^{23}
(q^{2\l})_3
X_V(\l,\mu)
(q^{-2\l})_3(\Cal R^{23})^{-1}.\tag 2.19\endgather
$$
The lemma is proved. 
$\square$\enddemo

\proclaim{Lemma 2.8}
$$
X_V(\l,\mu)=\Cal J^{3,12}(\l)
\Cal J^{12}(\l-\frac{1}{2} h^{(3)})
\Psi_V^{(1)}(\l+\frac{1}{2} (h^{(2)}-h^{(3)}),\mu)\Cal J^{32}(\l)^{-1}.
\tag 2.20
$$
\endproclaim

\demo{Proof}
Consider the function $X(\l)=\Cal J^{3,12}(\l)^{-1}X_V(\l)
\Cal J^{32}(\l)$. 
(for brevity we suppress the dependence on $\mu$ in the notation). 
By Lemma 2.7 and Lemma 2.4, this function satisfies the equation
$$
(q^{2\l})_3q^{\sum (x_i\o 1\o x_i+1\o x_i\o x_i)}X(\l)=X(\l)(q^{2\l})_3
q^{\sum 1\o x_i\o x_i}.\tag 2.21
$$
On the other hand, by the definition of $X_V$, the function 
$X$ is of the form \linebreak 
$Z_V(\l-\frac{1}{2} h^{(3)})+N$, where $N$ has negative
weights in the third component. As before, a solution of (2.21) 
with such property 
is unique, so it must coincide with its highest term. Thus, $N=0$,
and by Corollary 2.6 we get the Lemma. 
$\square$\enddemo

\proclaim{Corollary 2.9} 
$$
\gather
\Tr_0(\Cal R^{20}(\Cal R^{03})^{-1}\Phi^{V,01}_\mu q^{2\l} )=\\
(q^{-2\l})_2  
\Cal J^{3,12}(\l)
\Cal J^{12}(\l-\frac{1}{2} h^{(3)})
\Psi_V^{(1)}(\l+\frac{1}{2} (h^{(2)}-h^{(3)}),\mu)\Cal J^{32}(\l)^{-1}
(q^{2\l})_2.\tag 2.22
\endgather
$$
\endproclaim

\demo{Proof} Straightforward from Lemma 2.8.
$\square$\enddemo

\subhead 2.6. The difference operators\endsubhead

To obtain the action of the operators $\Cal M_W$ on 
traces, we need to apply to both sides of (2.22) 
the operation $m_{23}(1\o 1\o S)$, then multiply by $q^{2\rho}$ and take the 
trace.

Let $\Cal R=\sum a_i\o b_i$,
$\Cal R^{-1}=\sum a_i'\o b_i'$. 
Let $u=\sum S(b_i)a_i$ be the Drinfeld element. 
Let $\Cal J(\l)=\sum c_k\o d_k(\l)$, and
$\Cal J^{-1}(\l)=\sum c_k'\o d_k'(\l)$. Then we get
$$
\gather
\Cal M_W\Psi_V(\l,\mu)=\\
\Tr_2\sum [d_i^{(1)}(\l)c_j\o q^{-2\l}d_i^{(2)}(\l)d_j(\l+\frac{1}{2}h^{(2)})
]\Psi_V(\l+h^{(2)},\mu)(d_k'(\l)q^{2\l}S(c_k')S(c_i)q^{2\rho})_2.
\tag 2.23\endgather
$$
 
Let us simplify the expression for $\sum d_k'(\l)q^{2\l}S(c_k')$, 
which enters in (2.23).
 
\proclaim{Lemma 2.10} 
$$
\sum d_k'(\l)q^{2\l}S(c_k')=q^{\sum x_i^2}P(\l)S(u)q^{2\l},\tag 2.24
$$
where $P(\l):=\sum d_i'(\l)S^{-1}(c_i')$.
\endproclaim

\demo{Proof} From Lemma 2.4 we get 
$$
\Cal J^{-1}(\l)(q^{-2\l})_1=q^{-\sum x_i\o x_i}(q^{-2\l})_1\Cal J^{-1}(\l)
\Cal R^{21}.
$$
Applying $m_{21}(S\o 1)$ on both sides, we get
$$
\sum d_k'(\l)q^{2\l}S(c_k')=
q^{\sum x_i^2}\sum d_i'(\l)a_iS(b_i)S(c_i')q^{2\l},\tag 2.25
$$
which implies the Lemma, since $\sum a_iS(b_i)=S(u)$. 
$\square$\enddemo

Substituting (2.24) into (2.23) and using the cyclic property of the trace, 
we get
$$
\gather
\Cal M_W\Psi_V(\l,\mu)=\\
\sum d_i^{(1)}(\l)c_j\Psi_V(\l+\nu,\mu)
\Tr_2 [q^{\sum x_i^2}P(\l)S(u)q^{2\l}S(c_i)q^{2\rho}q^{-2\l}
d_i^{(2)}(\l)d_j(\l+\frac{1}{2}\nu)]=\\
\sum d_i^{(1)}(\l)c_j\Psi_V(\l+\nu,\mu)
\Tr_2[q^{\sum x_i^2}P(\l)q^{2\l}S^{-1}(c_i)q^{-2\l}
d_i^{(2)}(\l)d_j(\l+\frac{1}{2}\nu)S(u)q^{2\rho}]
\tag 2.26\endgather
$$

\proclaim{Lemma 2.11} 
$$
\gather
\sum d_i^{(1)}(\l)\o q^{2\l}S^{-1}(c_i)q^{-2\l}d_i^{(2)}(\l)=\\
\sum (a_j^{'(1)}d_k^{(1)}(\l))_1q^{-\sum (x_i)_2((x_i)_1+(x_i)_2)}
S^{-1}(c_k)S^{-1}(b_j')a_j^{'(2)}d_k^{(2)}(\l).\tag 2.27
\endgather
$$
\endproclaim

\demo{Proof} By Lemma 2.4, we have
$$
(q^{2\l})_1\Cal J^{1,23}(\l)(q^{-2\l})_1=(\Cal R^{23,1})^{-1}
\Cal J^{1,23}(\l)q^{\sum (x_i)_1((x_i)_2+(x_i)_3)}.\tag 2.28
$$
Applying $m_{13}(S^{-1}\o 1\o 1)$ to both sides, we get the lemma.
$\square$\enddemo

Now we will use the relations 
$$
\sum a_i^{'(1)}\o S^{-1}(b_i')a_i^{'(2)}=\sum a_k\o ub_k,
\tag 2.29
$$
$$
\Cal R^{23}\Cal J^{1,23}=\Cal J^{1,32}\Cal R^{23},\tag 2.30
$$
and 
$$
\sum S(c_i)d_i^{(1)}(\l)\o d_i^{(2)}(\l)=
S(\Q)(\l-\frac{1}{2}h^{(2)})_1\Cal J^{-1}(\l+\frac{1}{2}\hat h^{(1)}).\tag 2.31
$$
The first relation follows from the coproduct rule for the R-matrix, 
the second one is the functoriality of $\Cal J$, and the third one is obtained
if one applies $m_{12}(S\o 1\o 1)$ to the cocycle relation (2.12). 

Using (2.29)-(2.31), and the fact that $(1\o S)(\Cal R^{-1})=\Cal R$,
we get from (2.27) 
$$
\gather
\sum d_i^{(1)}(\l)\o q^{2\l}S^{-1}(c_i)q^{-2\l}d_i^{(2)}(\l)=\\
\sum (d_k^{(2)}(\l)a_j)_1(u^{-1})_2q^{-\sum (x_i)_2((x_i)_1+(x_i)_2)}
(S(c_k)d_k^{(1)}(\l)b_j)_2=\\
q^{-\sum (x_i)_2((x_i)_1+(x_i)_2)}(u^{-1})_2
S(\Q)(\l-\frac{1}{2}h)_2(\Cal J^{21})^{-1}(\l+\frac{1}{2}\hat h^{(2)})
\Cal R,
\tag 2.32\endgather
$$
where $\hat h^{(2)}$ should be replaced by the weight of the second 
component of the tensor product {\it after} the action of 
$\Cal J^{-1}$ (i.e. 
$$
\Cal J^{-1}(\l+\frac{1}{2}\hat h^{(2)})(v\o w)=
\sum_\nu \Cal J_\nu^{-1}(\l+\frac{1}{2}(wt(w)+\nu))(v\o w),
$$ 
where $wt(w)$ is the weight of $w$, and $\Cal J^{-1}_\nu$ is the part of 
$\Cal J^{-1}$ which shifts the weight of the second component by $\nu$). 

Substituting (2.32) into (2.26), and using the identity 
$u^{-1}S(u)=q^{-4\rho}$, we obtain 
$$
\Cal M_W\Psi_V(\l,\mu)=\sum_\nu\Tr|_{W[\nu]}(\tilde G(\l)\RR _{WV}(\l))
\Psi_V(\l+\nu,\mu),\tag 2.33
$$
where 
$$
\tilde G(\l)=q^{-2\rho}P(\l)S(\Q)(\l).\tag 2.34
$$
Thus Theorem 2.3 follows from the following Lemma.  

\proclaim{Lemma 2.12} $P(\l)=\Q^{-1}(\l+h)$. Hence 
$\tilde G(\l)=G(\l+h)$. 
\endproclaim

\demo{Proof} The Lemma is obtained by applying $m_{321}(S^{-1}\o 1\o S^{-1})$
to the cocycle identity for $\JJ $. 
$\square$\enddemo

Theorem 2.3 is proved. 

\subhead 2.7. Calculation of $G(\l)$\endsubhead

\proclaim{Proposition 2.13} Let $G(\l)=q^{-2\rho}\Q^{-1}(\l)S(\Q)(\l-h)$.  
Then 
$$
G(\l)=\frac{\delta_q(\l)}{\delta_q(\l-h)}.\tag 2.35
$$
\endproclaim

The rest of the subsection is the proof of Proposition 2.13.

\proclaim{Lemma 2.14}
$$
\Delta(G(\l))=\JJ (\l)(G(\l+h^{(2)})\o G(\l))\JJ ^{-1}(\l).\tag 2.36
$$
\endproclaim

\demo{Proof} The element $S(\Q)(\l)$ coincides with 
$K'(-\l-\rho)$, where the element $K'(\l)$
defined in Section 4.2 of \cite{EV}. By Lemma 27 of \cite{EV}, one has 
$\<w,K'(\l)w^*\>=B_{\l,W}(w,w^*)$, where the bilinear form 
$B_{\l,W}:W\o { ^*}W\to \C$ of weight zero is defined by the property 
$(1\o (,))\circ \Phi_{\l-\l_{w^*}}^w\Phi_\l^{w^*}=
B_{\l,W}(w,w^*)Id_{M_\l}$. Recall the 
functoriality property of $B_{\l,W}$:
$$
B_{U\o W}\circ (J_{UW}\o J_{{ ^*}W{ ^*}U})=B_U\circ B_W.\tag 2.37
$$
(for suitable $\l$-parameters, which are omitted for brevity).
This property implies that 
$$
\Delta(\Q(\l))=(S\o S)(\JJ ^{21}(\l)^{-1})
(\Q(\l)\o \Q(\l-h^{(1)}))\JJ(\l-h^{(1)}-h^{(2)})^{-1},
\tag 2.38
$$
which implies the lemma. 
$\square$\enddemo

Set $q=e^\gamma$, Then $J(\frac{\l}{\gamma})$ is a power series 
in $\gamma$ of the form $1+O(\gamma)$, whose terms lie in $U_q(\g)\o U_q(\g)
\o F$, where $F$ is the field of trigonometric functions.  

\proclaim{Lemma 2.15} Let $X(\l,\gamma)=\sum_{m=0}^\infty X_m(\l)\gamma^m$
be a power series of weight zero satisfying the equation 
$$
\Delta(X(\l,\gamma))=\JJ (\l/\gamma)(X(\l+\gamma h^{(1)},\gamma)\o 
X(\l,\gamma))
\JJ ^{-1}(\l/\gamma),\tag 2.39
$$
where $X_i\in U_q(\g)\o F$, 
such that $X_0(\l)=1$. Then:

(i) If $X=1$ on highest weight vectors
of all finite dimensional representations, then $X=1$. 

(ii) In general, $X=\frac{f(\l,\gamma)}{f(\l-\gamma h,\gamma)}$, where 
$f=1+\gamma f_1+\gamma^2 f_2+..$, $f_i\in F$.  
\endproclaim

\demo{Proof} (i) $X_1$ is a primitive element, so it lies in $\h\o F$.
 Since it acts trivially on all 
highest weight vectors, $X_1=0$. Similarly, $X_i=0$ for $i>1$ by 
induction, using the 
same argument. 

(ii) Let $\omega_i$ be fundamental weights, and $\eta_i(\l)$ be the 
eigenvalues of $X(\l)$ on the highest weight vector of the representation 
with highest weight $\omega_i$. It is easy to see that 
$\eta_i(\l+\gamma \omega_j)\eta_j(\l)=\eta_j(\l+\gamma\omega_i)\eta_i(\l)$, 
which implies that there exists $g$ such that 
$g(\l+\gamma\omega_i)=\eta_i(\l)g(\l)$. Let $\tilde X(\l)=
X(\l)g(\l)/g(\l+\gamma h)$. Then $\tilde X$ satisfies
the conditions of (i), so by (i) $\tilde X=1$. This proves (ii)
(we can set $f(\l)=g(\l+\gamma h)$). 
$\square$\enddemo

\proclaim{Corollary 2.16} 
$G(\l)=\frac{f(\l)}{f(\l-h)}$ for a suitable function $f$. 
\endproclaim

\demo{Proof} Follows from Lemmas 2.14 and 2.15. 
$\square$\enddemo

To conclude the proof of Proposition 2.13, it remains to show that
in Corollary 2.16, one can take $f(\l)=\delta_q(\l)$. 

To do this, we apply Theorem 2.3 and part (ii) of Proposition 2.2
in the case $V=\C$. 
In this case $\Psi_V(\l,\mu)=\frac{q^{2(\mu+\rho,\l)}}{\delta_q(\l)}$,
$R_{WV}=1$. So, using Corollary 2.16, we obtain
$$
\sum_\nu \frac{f(\l+\nu)}{f(\l)}|_{W[\nu]}\text{dim}W[\nu]
\frac{q^{2(\mu+\rho,\l+\nu)}}{\delta_q(\l+\nu)}=\chi_W(q^{2(\mu+\rho)})
\frac{q^{2(\mu+\rho,\l)}}{\delta_q(\l)}.\tag 2.40
$$
This equation is obviously satisfied if $f=\delta$. Since the validity of 
this equation for all $\mu$ completely determines $f(\l+h)/f(\l)$, 
we get the lemma.

\subhead 2.8. A modification of Theorem 2.3\endsubhead

Let $\varphi_V(\l,\mu)=\Psi_V(\l,\mu)\delta_q(\l)$.
 
\proclaim{Theorem 2.17} 
$$
\Cal D^{\l}_W\varphi_V(\l,\mu)=\chi_W(q^{2(\mu+\rho)})\varphi_V(\l,\mu).
\tag 2.41
$$
\endproclaim

\demo{Proof} Follows from Theorem 2.3, Proposition 2.2, and Proposition 2.13. 
$\square$\enddemo

\subhead 2.9. Proof of Theorem 1.1\endsubhead

The theorem follows from Theorem 2.17 and the fusion identity 
$$
\JJ ^{1...N}(\l)^{-1}\RR ^{0,1...N}(\l)\JJ ^{1...N}(\l+h^{(0)})=
\RR ^{01}(\l+h^{(2...N)})...\RR ^{0N}(\l),\tag 2.42
$$
which is easily checked from the definition.

\head 3. The dual MR equations.\endhead

\subhead 3.1. The map $L_{WV^*}(\mu)$\endsubhead

Let $W$ be a finite dimensional representation of $U_q(\g)$. 
Denote by $H_{\mu W}^{\l}$ the isotypic component 
of type $M_\l$ in $M_\mu\o W$. For generic $\mu$, we have a natural 
isomorphism $\eta: W[\l-\mu]\o M_\l\to H_{\mu W}^\l$ given by
$w\o v\to \Phi_\l^wv$.  

Consider the intertwining operator 
$$
P_{V\o V^*,W}\Cal R_{VW}(\Phi_\mu^V\o Id_W):M_\mu\o W
\to M_\mu\o W\o V\o V^*.\tag 3.1
$$
Restricting this operator to the isotypic component of $M_{\mu+\nu}$, 
from the intertwining property we obtain 
$$
P_{V\o V^*,W}\Cal R_{VW}(\Phi_\mu^V\o Id_W)|_{W[\nu]\o M_{\mu+\nu}}=
L_{WV^*}(\mu)(Id_W\o \Phi_{\mu+\nu}^V),\tag 3.2
$$
where $L_{WV^*}(\mu):W\o V^*\to W\o V^*$ is a uniquely determined 
linear map. Our task now is to find this linear map explicitly.

\subhead 3.2. Computation of $L_{WV^*}$ \endsubhead

\proclaim{Proposition 3.1}
$$
L_{WV^*}(\mu)=
R_{WV}(\mu+\nu)^{t_2},\tag 3.3
$$
where $t_2$ means transposition in the third component. 
\endproclaim

\demo{Proof} Let us apply both sides of (3.2) to $w\o y$. 
Let $B$ be a homogeneous basis of $V$. 
Then we get 
$$
\sum_{v\in B}[\eta^{-1}P_{VW}\Cal R_{VW}
(\Phi_\mu^v\o 1)\Phi_{\mu+\nu}^wy]\o v^*
=\sum_{v\in B}L_{WV^*}(\mu)\Phi_{\mu+\nu}^v(w\o y\o v^*).\tag 3.4
$$
Moving $\eta^{-1}$ and $\Cal R$ to the right, and using the notation 
$L_{WV^*}=\sum p_i\o q_i$, we get
$$
\sum_{v\in B}
[(\Phi_\mu^v\o 1)\Phi_{\mu+\nu}^wy]\o v^*
=P_{VW}(\Cal R_{VW}^{21})^{-1}\sum_{i}\sum_{v\in B}
\Phi_{\mu+\nu-wt(v)}^{p_iw}
\Phi_{\mu+\nu}^{v}y\o 
q_iv^*,
\tag 3.5
$$
where $wt(v)$ denotes the weight of $v$. 

Simplifying (3.5), we get
$$
\sum_{v\in B}
\Phi_{\mu+\nu}^{J(\mu+\nu)(v\o w)}y\o v^*
=\sum_{i}\sum_{v\in B}
\Phi_{\mu+\nu}^{P_{VW}(\Cal R^{21})^{-1}
J(\mu+\nu)(p_iw\o q_i^*v)}y\o 
v^*.
\tag 3.6
$$
This implies that on $W[\nu]\o V$, one has
$$
J(\mu+\nu)=\Cal R^{-1}J^{21}(\mu+\nu)(L^{21}_{WV^*}(\mu))^{t_1}.\tag 3.7
$$
That is, 
$$
L_{WV^*}(\mu)=(J^{-1}(\mu+\nu)\Cal R^{21}J^{21}(\mu+\nu))^{t_2}=
R_{WV}(\mu+\nu)^{t_2}.\tag 3.8
$$
$\square$\enddemo

\subhead 3.3. Proof of Theorem 1.2\endsubhead

Let us multiply both sides of (3.2) on the right by $q^{2\l}$ (acting 
in $M_\mu\o W$), on the left by 
$q^{-\sum x_i|_V x_i|_W}$, and compute the trace in the Verma modules, 
summing over
all $\nu$. After multiplication by $\delta_q(\l)$
and using 
Proposition 3.1 and the fact that $\varphi_V$ 
has zero weight in $V$, we obtain
$$
\chi_W(q^{2\l})\varphi_V(\l,\mu)=\sum_\nu \Tr|_{W(\nu)}
(R_{WV}(\mu+\nu)^{t_2})\varphi_V(\l,\mu+\nu).\tag 3.9
$$
This is equivalent to 
$$
\chi_{W^*}(q^{-2\l})\varphi_V(\l,\mu)=\sum_\nu \Tr|_{W^*(-\nu)}
(R_{WV}(\mu+\nu)^{t_1t_2})\varphi_V(\l,\mu+\nu).\tag 3.10
$$
Rewriting the last equation in terms of $F_V(\l,\mu)$, we obtain
$$
\chi_{W^*}(q^{-2\l})F_V(\l,\mu)=
\sum_\nu \Tr|_{W^*(\nu)} \Q^{-1}(\mu)|_{V^*}
(\RR_{WV}(\mu+\nu)^{t_1t_2})\Q(\mu+\nu)|_{V^*}F_V(\l,\mu+\nu).\tag 3.11
$$
The expression $\RR_{WV}(\mu)^{t_1t_2}$ can be computed from (2.38):
$$
\RR_{WV}(\mu)^{t_1t_2}=(\Q(\mu)\o \Q(\mu-h^{(1)}))
\RR_{W^*V^*}(\mu-h^{(1)}-h^{(2)})(\Q^{-1}(\mu-h^{(2)})\o \Q^{-1}(\mu)).\tag 3.12
$$
This expression and (3.11) imply Theorem 1.2.

\head 4. The dual qKZB equations \endhead 

In this section, we will prove Theorem 1.4. 

 Using formula (3) in \cite{EV}, we get
$$
(\Phi_{\mu-\mu_v}^{w}\o 1)\Phi_{\mu}^{v}=
(1\o \Cal R_{WV}^{-1}P_{VW})\sum_k 
(\Phi_{\mu-\mu_{w_k}}^{v_k}\o 1)\Phi_{\mu}^{w_k},\tag 4.1
$$
where 
$$
\sum v_k\o w_k=R_{VW}(\mu)(v\o w).\tag 4.2
$$ 
Therefore, using simplified notations, we have
$$
\Phi_{\mu+h^{(V^*)}}^{W}\Phi_{\mu}^{V}=\Cal R^{-1}
R_{21}^*(\mu)\Phi_{\mu+h^{(W^*)}}^{V}\Phi_{\mu}^{W}=
\Cal R_{21}
R^*(\mu)^{-1}\Phi_{\mu+h^{(W^*)}}^{V}\Phi_{\mu}^{W}
\tag 4.3
$$

Let us now 
take the j-th intertwining operator in (1.3) and move it to the right, 
permuting it with other operators. 
We have 
$$
\gather
\Psi_{V_1...V_N}(\l,\mu)=\Cal R_{j+1j}...
\Cal R_{Nj}(q^{2\l})_j
R_{jj+1}^*(\mu+\sum_{i=j+2}^Nh^{(*i)})^{-1}...R_{jN}^*(\mu)^{-1}
\times \\
\text{Tr}(\Phi^{V_1}_{\mu+\sum_{i=2}^Nh^{(*i)}}
...\Phi^{V_{j-1}}_{\mu+\sum_{i=j}^Nh^{(*i)}}
\Phi^{V_{j+1}}_{\mu+h^{(*j)}+\sum_{i=j+2}^Nh^{(*i)}}
...\Phi^{V_N}_{\mu+h^{(*j)}} q^{2\l}\Phi^{V_j}_\mu)=\\
\Cal R_{j+1j}...
\Cal R_{Nj}(q^{2\l})_j
R_{jj+1}^*(\mu+\sum_{i=j+2}^Nh^{(*i)})^{-1}...R_{jN}^*(\mu)^{-1}
\Gamma_{*j}
\times \\
\text{Tr}(\Phi^{V_1}_{\mu+\sum_{i=2,i\ne j}^Nh^{(*i)}}
...\Phi^{V_{j-1}}_{\mu+\sum_{i=j+1}^Nh^{(*i)}}
\Phi^{V_{j+1}}_{\mu+\sum_{i=j+2}^Nh^{(*i)}}
...\Phi^{V_N}_{\mu} q^{2\l}\Phi^{V_j}_{\mu+\sum_{i\ne j}h^{(*i)}})
\tag 4.4
\endgather
$$
(in the last equality we use that $\sum h^{(*i)}=0$). 
Using the cyclic property of the trace, we can now 
put the j-th operator into the 
beginning, and move it to the right to its original place, thus completing 
the cycle. This yields
$$
\gather
\Psi_{V_1...V_N}(\l,\mu)=
\Cal R_{j+1j}...\Cal R_{Nj}(q^{2\l})_j\Cal R_{j1}^{-1}...
\Cal R_{jj-1}^{-1}\times \\
R_{jj+1}^*(\mu+\sum_{i=j+2}^N
h^{(*i)})^{-1}...
 R_{jN}^*(\mu)^{-1}\Gamma_{*j}R_{1j}^*(\mu+\sum_{i=2,i\ne j}^Nh^{(*i)})...
R_{j-1j}^*(\mu+\sum_{i=j+1}^Nh^{(*i)})
\times \\
\Psi_{V_1...V_N}(\l,\mu)
.\tag 4.5
\endgather
$$

Multiplying both sides of (4.5) by $\delta_q(\l)$, 
replacing $\mu$ by $-\mu-\rho$, and multiplying 
by $\JJ^{1...N}(\l)^{-1}$ and the product of $\Q^{-1}$-s,
we obtain
$$
\gather
F_{V_1...V_N}(\l,\mu)=
\JJ^{1...N}(\l)^{-1}
\Cal R_{j+1j}...\Cal R_{Nj}(q^{2\l})_j\Cal R_{j1}^{-1}...
\Cal R_{jj-1}^{-1} 
\JJ^{1...N}(\l)\\
[\Q^{-1}(\mu)^{(*N)}\o ...\o \Q^{-1}(\mu-h^{(*2...*N)})^{(*1)}]\times \\
\RR_{jj+1}^*(\mu-\sum_{i=j+2}^N
h^{(*i)})^{-1}...
 \RR_{jN}^*(\mu)^{-1}\Gamma_{*j}^{-1}
\RR_{1j}^*(\mu-\sum_{i=2,i\ne j}^Nh^{(*i)})...
\RR_{j-1j}^*(\mu-\sum_{i=j+1}^Nh^{(*i)})
\times \\
[\Q(\mu)^{(*N)}\o ...\o \Q(\mu-h^{(*2...*N)})^{(*1)}]
F_{V_1...V_N}(\l,\mu)
.\tag 4.6
\endgather
$$
 
Inverting (4.6), we get
$$
\gather
F_{V_1...V_N}(\l,\mu)=
\JJ^{1...N}(\l)^{-1}
\Cal R_{jj-1}...\Cal R_{j1}(q^{-2\l})_j\Cal R_{Nj}^{-1}...
\Cal R_{j+1j}^{-1} 
\JJ^{1...N}(\l)\\
[\Q^{-1}(\mu)^{(*N)}\o ...\o \Q^{-1}(\mu-h^{(*2...*N)})^{(*1)}]\times \\
\RR_{j-1j}^*(\mu-\sum_{i=j+1}^Nh^{(*i)})^{-1}...
 \RR_{1j}^*(\mu-\sum_{i=2,i\ne j}^Nh^{(*i)})^{-1}
\Gamma_{*j}\RR_{jN}^*(\mu)...
\RR_{jj+1}^*(\mu-\sum_{i=j+2}^Nh^{(*i)})
\times \\
[\Q(\mu)^{(*N)}\o ...\o \Q(\mu-h^{(*2...*N)})^{(*1)}]
F_{V_1...V_N}(\l,\mu)
.\tag 4.7
\endgather
$$

Using formula (3.12), it is not difficult to check that 
on zero weight vectors
$$
\gather
[\Q^{-1}(\mu)^{(*N)}\o ...\o \Q^{-1}(\mu-h^{(*2...*N)})^{(*1)}]\times \\
\RR_{j-1j}^*(\mu-\sum_{i=j+1}^Nh^{(*i)})^{-1}...
 \RR_{1j}^*(\mu-\sum_{i=2,i\ne j}^Nh^{(*i)})^{-1}
\Gamma_{*j}\RR_{jN}^*(\mu)...
\RR_{jj+1}^*(\mu-\sum_{i=j+2}^Nh^{(*i)})
\times \\
[\Q(\mu)^{(*N)}\o ...\o \Q(\mu-h^{(*2...*N)})^{(*1)}]=K_j^\vee.\tag 4.8
\endgather
$$

So it remains to check that 
$$
\JJ^{1...N}(\l)^{-1}
\Cal R_{jj-1}...\Cal R_{j1}(q^{-2\l})_j\Cal R_{Nj}^{-1}...
\Cal R_{j+1j}^{-1} 
\JJ^{1...N}(\l)=D_j^\vee\tag 4.9
$$
(on zero weight vectors).

To check (4.9), observe that it can be written in the form
$$
\JJ^{1...N}(\l)^{-1}
\Cal R_{j,1...j-1}(q^{-2\l})_j\Cal R_{j+1...N,j}^{-1} 
\JJ^{1...N}(\l)=D_j^\vee,\tag 4.10
$$
which implies that it is enough to check (4.9) for $N=2$ and $N=3$. 

Let $E_j$ be the left hand side of (4.9). Since 
$E_1E_2E_3=D_1^\vee D_2^\vee D_3^\vee =1$ for $N=3$, it is enough to prove 
(4.9) for $N=2$ and $N=3$, $j=1,3$. But the last two cases follow from the 
case $N=2$, so it is enough to check only this case. 

If $N=2$, identity (4.9) easily follows from Lemma 2.4. 
Theorem 1.4 is proved. 

\head 5. The symmetry identity.\endhead

In this section we will prove Theorem 1.5 using Theorems 1.1 and 1.2. 

Let $V$ be a finite dimensional representation of $U_q(\g)$, and 
$v\in V$ a homogeneous vector. 
Consider the operator $\Phi_\mu^v$ as a linear operator 
$U_q(\n_-)\to U_q(\n_-)\o V$, identifying the Verma modules
with $U_q(\n_-)$. 

\proclaim{Lemma 5.1}  The operator 
$\Phi_\mu^v$ is a rational function of variables $q^{2(\mu,\alpha_i)}$
of the form $\tilde\Phi_\mu^v/D(\mu)$, where 
$\tilde\Phi_\mu^V$ is an operator valued 
polynomial in $x_i:=q^{2(\mu,\alpha_i)}$, and 
$D(\mu)$ is a polynomial in $x_i$ with nonzero free term. 
\endproclaim

\demo{Proof}
This follows from the arguments in the proof of Proposition 2.2 in [ES].
$\square$\enddemo

\proclaim{Corollary 5.2} The matrix elements of the universal 
trace function $\Psi$ can be expanded in formal series 
belonging to the space $q^{2(\l,\mu)}\C[[q^{-2(\l,\alpha_i)}]]\o
\C[[q^{2(\mu,\alpha_i)}]]$.
\endproclaim

{\bf Remark.} Note that the tensor product in Corollary 
5.2 is algebraic, i.e. uncompleted.

Let $K=q^{-2(\l,\mu)}\C[[q^{-2(\l,\alpha_i)}]]\o
\C[[q^{-2(\mu,\alpha_i)}]]$.
Denote by $L$ the space \linebreak $(V_1\o...\o V_N)[0]\o K$. 

\proclaim{Lemma 5.3} 
For any vector $v\in V_1\o...\o V_N[0]$, 
the function $\<F_{V_1...V_n},v\>$ belongs to $L$. 
\endproclaim 

\demo{Proof}
It is clear from Proposition 2.2 of \cite{ES} 
and the definition of $\JJ$ that for any $V,W$ the function $\JJ_{VW}(\l)$ 
is a rational function of $y_i=q^{-2(\l,\alpha_i)}$, 
which can be represented as a ratio of an operator-valued polynomial 
of $y_i$ whose free term is invertible,
and a scalar polynomial of $y_i$ with nonzero free term. The same statement, 
for the same reason, applies to $\Q(\l)$. These 
facts and Corollary 5.2 imply the Lemma.
$\square$\enddemo

Let $L_0\subset L$ be the space of solutions 
of equations (1.11) whose highest term (as of a series in $x_i=
q^{-2(\l,\alpha_i)}$) is $v q^{-2(\l,\mu)}$, where $v\in (V_1\o...V_N)[0]$
is independent of $\mu$.  
By Theorem 1.1 and Lemma 5.3, 
for any $v\in (V_1\o...\o V_N)[0]$ the inner product 
$\<F_{V_1,...,V_N},v\>$ is in $L_0$. 

\proclaim{Lemma 5.4} 
$L_0$ is a complex vector space 
of dimension $d=\text{dim}(V_1\o...\o V_N)[0]$. 
Moreover, if $B$ is a basis of $(V_1\o...\o V_N)[0]$ then the collection of 
functions $\<F_{V_1,...,V_N},v\>$, $v\in B$, is a basis of $L_0$. 
\endproclaim

\demo{Proof} The solutions $\<F_{V_1,...,V_N},v\>,v\in B$ are linearly 
independent, since, by the arguments in the proof of Lemma 5.3, 
 so are their highest terms. So it remains to show 
that the dimension of the space $L_0$ is not bigger than $d$. 

Let $I$ be the maximal ideal in $\C[[x_1,...,x_r]]$. It is clear that 
for any $m\ge 1$ the operator $\Cal D_W-\chi_W(q^{-2\mu})$ 
preserves the subspace
$I^mL$, so it is enough to check that the dimension 
of the space of solutions is at most $d$ on $L/I^mL$ for all $m$. 

For $s\in \C^*$, define an automorphism $g_s$ of $L$ by 
$g_sf(x_1,...,x_r)=f(sx_1,...,sx_r)$. Let 
$\Cal D_W^0=\lim_{s\to 0}g_s\Cal D_Wg_s^{-1}$. 
It is enough to show that the system 
$\Cal D_W^0f=\chi_W(q^{-2\mu})f$ has no more than $d$ linearly 
independent solutions. 

The operator $D_W^0$ is easily computed. Namely, since
the exchange matrices converge to the usual R-matrices
as $\l\to\infty$
(see Theorem 50 of \cite{EV}), one gets
$$
\Cal D_W^0=\sum_\nu \text{dim}(W[\nu])T_\nu.
$$
Therefore, the operators $\Cal D_W^0$ 
are diagonal in the monomial basis, and 
it is obvious
that the system 
$\Cal D_W^0f=\chi_W(q^{-2\mu})f$ has no more than $d$ linearly 
independent solutions.
$\square$\enddemo

\proclaim{Corollary 5.5} There exists a unique element
$M(\mu)\in End((V_N^*\o...\o V_1^*)[0])\o 
\C[[q^{-2(\mu,\alpha_i)}]]$ such that 
$$
F^*_{V_N^*,...,V_1^*}(\mu,\l)=(1\o M(\mu))F_{V_1,...,V_N}(\l,\mu).\tag 5.2
$$
\endproclaim

\demo{Proof} This follows from Lemma 5.4 and Theorem 1.2.
$\square$\enddemo

Now we will prove Theorem 1.5. Applying Corollary 5.5 twice, we 
obtain that $M(\mu)M'(\l)=1$ for some function $M'(\l)$. This implies that $M$
is in fact independent of $\mu$. 
Looking at the highest terms of $F$ and $F^*$, one finds that $M=1$.  

This proves Theorem 1.5. 

\head 6. The symmetry of the universal trace function under 
$q\to q^{-1}$, and the function $u_V$. 
\endhead

\subhead 6.1. The symmetry identity with $q\to q^{-1}$\endsubhead

In this section we will study the behavior of the functions 
$\Psi_V$
under the transformation $q\to q^{-1}$. 

Consider the algebra isomorphism $\xi:U_q(\g)\to U_{q^{-1}}(\g)$ 
given by $\xi(E_i)=E_iq^{d_ih_i}$, $\xi(F_i)=q^{-d_ih_i}F_i$, 
$\xi (q^h)=q^{-h}$. Using this isomorphism, we will identify 
$U_q(\g)$ and $U_{q^{-1}}(\g)$ (as algebras), and regard 
any module over $U_q(\g)$ as a module over $U_{q^{-1}}(\g)$ and vice versa.

We will start by proving the following theorem. 

\proclaim{Theorem 6.1} The function $\Psi_{V_1,..,V_N}$ satisfies 
the equality
$$
\Psi_V(q^{-1},-\l,\mu)=
u(q)^{-1}\Q(q,\l)\Psi_V(q,\l,\mu),
\tag 6.1
$$
where $u(q)$ denotes the Drinfeld element 
$u$ of $U_{q}(\g)$. 
\endproclaim

{\bf Remark.} To avoid confusion, here and below we 
explicitly specify for every expression 
whether it is evaluated at $q$ or at $q^{-1}$.

\demo{Proof} 
It is obvious that $\xi$ is an antiisomorphism of 
coalgebras. This implies that for a finite dimensional representation 
$V$ and a vector $v\in V[0]$ the operator 
$P\Phi_\mu^v(q)$ is an intertwining operator $M_\mu\to V\o M_\mu$
for $U_{q^{-1}}(\g)$, and hence
the operator $(\Cal R^{10})^{-1}\Phi_\mu^v(q)$ 
is an intertwining operator $M_\mu\to M_\mu\o V$ for $U_{q^{-1}}(\g)$. 
We have $(\Cal R^{10})^{-1}(q)\Phi_\mu^v(q) v_\mu=v_\mu\o v+...$, so 
we have 
$$
(\Cal R^{10})^{-1}(q)\Phi_\mu^V(q)=\Phi_\mu^V(q^{-1}).\tag 6.2
$$
Let us multiply both sides of (6.2) by $q^{2\l}$ and take the trace. 

Using (2.15) and (2.16), we get 
$$
\text{Tr}((\Cal R^{10})^{-1}(q)\Phi_\mu^V(q)q^{2\l})=
u(q)^{-1}\Q(q,\l)\Psi_V(q,\l,\mu),\tag 6.3
$$
which implies the theorem. 
$\square$\enddemo

\subhead 6.2. The function $\hat u_V(\l,\mu)$\endsubhead 

Now let $V$ be an irreducible $U_q(\g)$-module, and define  
the function 
$$
\hat u_V(\l,\mu):=\delta_q(\l)\Psi_V(q^{-1},-\l,-\mu-\rho).\tag 6.4
$$

\proclaim{Corollary 6.2}
The function $\hat u_V$ is symmetric: 
$\hat u_V(\l,\mu)=\hat u_{V^*}^*(\mu,\l)$. 
\endproclaim

\demo{Proof} Let $\nu$ be the highest weight of $V$. 
By Theorem 6.1, we have 
$$
\hat u_V(\l,\mu)=q^{-(\nu,\nu+2\rho)}\delta_q(\l)\Q(q,\l)\Psi_V(\l,-\mu-\rho)=
q^{-(\nu,\nu+2\rho)}(\Q(q,\l)\o \Q(q,\mu))F_V(\l,\mu).\tag 6.5
$$
Therefore, the symmetry of $\hat u_V$ follows from Theorem 1.5. 
$\square$\enddemo

\subhead 6.3. The function $u_V(\l,\mu)$\endsubhead

The function $\hat u_V$ has poles. We would like to get rid
of them by multiplying the function $\hat u_V$ by its denominator. 

Following \cite{ES}, introduce the function
$$
\delta_V(\mu)=\prod_{\alpha>0}\prod_{n=1}^{k_\alpha^V}
(q^{(\alpha,\mu+\rho)-n(\alpha,\alpha)/2}-
q^{-(\alpha,\mu+\rho)+n(\alpha,\alpha)/2}),\tag 6.6
$$
where for a positive root $\alpha$ we define
$k_\alpha^V:=\max 
\{n: V[n\alpha]\ne 0\}$.

Define the function 
$$
u_V(\l,\mu):=\delta_{V^*}(-\l-\rho)\delta_V(-\mu-\rho)
\hat u_V(\l,\mu).\tag 6.7
$$

\proclaim{Proposition 6.3} 
The function $u_V$ is symmetric: 
$u_V(\l,\mu)=u_{V^*}^*(\mu,\l)$, and is a product of 
$q^{-2(\l,\mu)}$ and a trigonometric polynomial of $\l$ and $\mu$. 
In particular, it is holomorphic in $\l,\mu$.  
\endproclaim

\demo{Proof} The symmetry of $u_V$ follows from Corollary 6.2.
The fact that $u_V$ is holomorphic in $\mu$ is a consequence of 
Lemma 5.1. The fact that $u_V$ is holomorphic in $\l$ 
follows from the symmetry. 
$\square$\enddemo

Using (6.5), one gets the following expression of $u_V$ in terms of $F_V$:
$$
u_V(\l,\mu)=q^{-(\nu,\nu+2\rho)}(\delta_{V^*}(-\l-\rho)\Q(q,\l)\o 
\delta_V(-\mu-\rho)\Q(q,\mu))F_V(\l,\mu).\tag 6.8
$$

In Section 8 we will show that the function $u_V$ for $U_q(sl_2)$ coincides 
(up to a constant factor) with 
the trigonometric limit of the universal 
hypergeometric function $u$ introduced in 
\cite{FV2}.

\head 7. Calculation of the functions $u_V$ and $F_V$ for $sl_2$\endhead

\subhead 7.1. Calculation of the trace function\endsubhead

Recall that $U_q(sl_2)$ is generated by $E,F,q^h$, with relations 
$$
q^hEq^{-h}=q^2E,\ q^hFq^{-h}=q^{-2}F,\ EF-FE=\frac{q^h-q^{-h}}{q-q^{-1}},
\tag 7.1
$$
and the coproduct is defined by
$$
\Delta(E)=E\o q^h+1\o E,\ \Delta(F)=F\o 1+q^{-h}\o F. \tag 7.2
$$

Weights for $U_q(sl_2)$ can be identified with complex numbers: 
we say that a vector $v$ in an $U_q(sl_2)$-module 
has weight $\mu$ if $q^hv=q^\mu v$. In this case we have 
$(\mu,\mu')=\mu\mu'/2$.
We will write $q^\l$ for $q^{(\l,\alpha)}$. 
Thus the meaning of $q^\l$ in this section is different from the 
previous sections.  

Recall also that for any number $a$ the q-number $[a]_q$ is defined by 
$[a]_q=\frac{q^a-q^{-a}}{q-q^{-1}}$. 

Consider the function $\Psi_V(\l,\mu)$, where $V$ is 
the representation of $U_q(sl_2)$ with highest weight $2m$, $m\in\Z_+$. 
Since the weight space $V[0]$ is 1-dimensional, 
this function can be considered as a scalar function. 

\proclaim{Theorem 7.1} The function $\Psi_V(\l,\mu)$ is given by 
the formula
$$
\gather
\Psi_V(\l,\mu)=\\
q^{\l\mu}
\sum_{l=0}^mq^{l(l-1)/2}(q-q^{-1})^l 
\frac{[m+l]_q!}{[l]_q![m-l]_q!}
\frac{q^{-2l\l}}{\prod_{j=0}^{l-1}(1-q^{2(\mu-j)})
\prod_{j=0}^l(1-q^{-2(\l-j)})},\tag 7.3
\endgather
$$
where $[n]_q!=[1]_q...[n]_q$. 
\endproclaim

\demo{Proof} We fix a generator $w_0$ of $V[0]$. 
Let us compute the intertwining 
operator $\Phi_\mu:M_\mu\to M_\mu\o V$.

Let $w_\beta, \beta=m,m-1,...,-m$, be the basis of $V$ defined by 
the condition $Fw_\beta=w_{\beta-1}$ if $\beta\ne -m$
(this basis is unique up to a common scalar). 
The image of the highest weight vector under 
the operator $\Phi_\mu$ has the form 
$$
\Phi_\mu v_\mu=\sum_{j=0}^mc_j(\mu)F^jv_\mu\o w_j,\tag 7.4
$$
where $c_0=1$.

Let us compute the coefficients $c_j(\mu)$. They are computed from the 
condition that $\Delta(E)=E\o q^h+1\o E$ annihilates the r.h.s.
of (7.4). Using the formula
$$
ew_\beta=[m+\beta+1]_q[m-\beta]_qw_{\beta+1},\ \beta\ne m,\tag 7.5
$$
we can rewrite this condition in the form
$$
q^{2i}[i]_q[\mu-i+1]_qc_i(\mu)+[m+i]_q[m-i+1]_qc_{i-1}(\mu)=0,\ i\ge 1.\tag 7.6
$$
 Solving this recurrence relation, we obtain
$$
c_i(\mu)=(-1)^iq^{-i(i+1)}\frac{[m+i]_q!}{[i]_q![m-i]_q!}
\prod_{j=1}^i[\mu-i+1]_q^{-1}.\tag 7.7
$$

Now we need to compute $\Phi_\mu F^kv_\mu$. For this, we need to compute 
$\Delta(F^k)$. Using the q-binomial theorem, we obtain
$$
\Delta(F^k)=(F\o 1+q^{-h}\o F)^k=
\sum_{l=0}^k \left( \matrix k\\ l\endmatrix\right)_{q^{-2}}
q^{-lh}F^{k-l}\o F^l,\tag 7.8
$$
where 
$$
\left(\matrix k\\ l\endmatrix\right)_{p}:=
\frac{\prod_{i=k-l+1}^{k}(1-p^{i})}{\prod_{i=1}^{l}(1-p^{i})}.\tag 7.9
$$

Now, using the intertwining property of 
$\Phi_\mu$, we obtain
$$
\gather
\Phi_\mu F^kv_\mu=\Delta(F^k)\Phi_\mu v_\mu=\\
(\sum_{l=0}^k \left(\matrix k\\ l \endmatrix\right)_{q^{-2}}
q^{-lh}F^{k-l}\o F^l)\sum_{j=0}^mc_j(\mu)F^jv_\mu\o w_j
,\tag 7.10
\endgather
$$
This double sum reduces to a single sum if we use 
the following version of the q-binomial theorem:
$$
\sum_{k\ge l} 
\left(\matrix k\\ l\endmatrix\right)_{p} x^{k-l}
=\prod_{i=0}^l(1-p^ix)^{-1}.\tag 7.11
$$
Substituting (7.7) into (7.10) and using (7.11),
we obtain (7.3). The Theorem is proved. 
$\square$\enddemo

\proclaim{Corollary 7.2} 
$$
\gather
u_V(\l,\mu)=\\
q^{-\l\mu}
\sum_{l=0}^mq^{-l(l-1)/2}(q-q^{-1})^l 
\frac{[m+l]_q!}{[l]_q![m-l]_q!}
q^{-l(\l+\mu)}\prod_{j=l+1}^{m}(q^{\l+j}-q^{-\l-j})(q^{\mu+j}-q^{-\mu-j}).
\tag 7.12\endgather
$$
\endproclaim

\demo{Proof} The statement is obtained from Theorem 7.1 and formulas
(6.4) and (6.7).
$\square$\enddemo

The function $u_V$ is manifestly symmetric in $\l$ and $\mu$, as predicted by 
Proposition 6.3.  

\subhead 7.2. A formula for $\Q(\mu)$\endsubhead

Now we would like to compute the function $F_V$. So it
remains to compute the value 
of $\Q(\mu)$ on the zero-weight subspace of $V$. 
This value is given by the following Lemma.

\proclaim{Lemma 7.2} 
The element $\Q(\mu)$ acts on the zero weight subspace of $V$ by 
the formula
$$
\Q(\mu)|_{V[0]}=
q^{-2m}\prod_{j=1}^{m}
\frac{q^{-2\mu-2j+2}-q^{-2m}}{q^{-2\mu-2j}-1}.\tag 7.13
$$

\endproclaim

\demo{Proof} Denote by $\Q_{r,l}(\mu)$ the eigenvalue of $\Q(\mu)$ on the
subspace of weight $l$ of the representation of $U_q(sl_2)$ with 
highest weight $r$. 

It is clear that $\Q_{0,0}=1$. The values of $\Q_{1,1}$ and 
$\Q_{1,-1}$ are easily computed from the definition. Namely, 
from the ABRR equation we have 
$$
\Cal J(\mu)=1+\frac{q^{-1}-q}{q^{-2\mu}-q^{-2}(q^{-h}\o q^h)}F\o E+...,
$$
and therefore
$$
\Q_{1,-1}=1,\ \Q_{1,1}=\frac{q^{-2\mu}-q^{-2}}{q^{-2\mu}-1}.\tag 7.14
$$

Now consider the subspace of weight $l$ in the tensor product 
$\C^2\o W$, where $W$ has highest weight $r$. 
Take the determinant of both sides of (2.38) 
restricted to this subspace. Using the strict triangularity of 
$\JJ$, we can ignore the $\JJ$ terms and obtain
$$
\Q_{r+1,l}(\mu)\Q_{r-1,l}(\mu)=\Q_{r,l-1}(\mu-1)\Q_{r,l+1}(\mu+1) 
\frac{q^{-2\mu}-q^{-2}}{q^{-2\mu}-1}.\tag 7.15
$$
This implies 
$$
\frac{\Q_{r+1,l}(\mu)}{\Q_{r,l+1}(\mu+1)}=
\frac{\Q_{r,l-1}(\mu-1)}{\Q_{r-1,l}(\mu)} 
\frac{q^{-2\mu}-q^{-2}}{q^{-2\mu}-1},\tag 7.16
$$
from which we get
$$
\frac{\Q_{r+1,l}(\mu)}{\Q_{r,l+1}(\mu+1)}=
\prod_{j=0}^{(r+l-1)/2}
\frac{q^{-2\mu+2j}-q^{-2}}{q^{-2\mu+2j}-1}=
\frac{q^{-2\mu}-q^{-r-l-1}}{q^{-2\mu}-1}
\tag 7.17
$$
This yields
$$
\Q_{r+1,l}(\mu)=\prod_{j=0}^{(r-l+1)/2}
\frac{q^{-2\mu-2j}-q^{-r-l-1}}{q^{-2\mu-2j}-1}.\tag 7.18
$$
In particular, 
$$
\Q_{2m,0}(\mu)=\prod_{j=0}^{m}
\frac{q^{-2\mu-2j}-q^{-2m}}{q^{-2\mu-2j}-1}=
q^{-2m}\prod_{j=1}^{m}
\frac{q^{-2\mu-2j+2}-q^{-2m}}{q^{-2\mu-2j}-1},\tag 7.19
$$
as desired.
$\square$\enddemo

\subhead 7.3. Calculation of 
$F_V(\l,\mu)$\endsubhead

\proclaim{Proposition 7.3}
The function $F_V(\l,\mu)$ is given by 
$$
\gather
F_V(\l,\mu)=q^{-\l\mu}\prod_{j=1}^{m}
\frac{q^{-2\mu-2j}-1}{q^{-2\mu-2j+2}-q^{-2m}}\times\\
q^{2m}\sum_{l=0}^mq^{l(l-1)/2}(q-q^{-1})^l 
\frac{[m+l]_q!}{[l]_q![m-l]_q!}
\frac{q^{-2l\l}}{\prod_{j=1}^{l}(1-q^{-2(\mu+j)})
\prod_{j=1}^l(1-q^{-2(\l-j)})},
\tag 7.20
\endgather
$$
\endproclaim

\demo{Proof}
Using the definition of $F_V$, Theorem 7.1 and Lemma 7.2, we get
(7.20).
$\square$\enddemo

{\bf Remark.} Note that expression (7.20) is not manifestly symmetric. 
In fact, the separate terms in the sum (7.20) are not 
symmetric, and it is only the whole sum that has the symmetry $\l\to\mu$.

\head 8. Integral representation of the trace function $u_V$ 
for $\g=sl_2$\endhead

In \cite{FV2}, G.Felder and the second author, studying the qKZB difference 
equations, defined 
the universal hypergeometric function $u_m(\l,\tau,\mu,p)$ 
(depending on a parameter $q$) with a number of interesting 
properties. In this section we will consider the trigonometric limit 
of this function, and will show that it coincides, up to a constant factor, 
with the 
function $u_{V(m)}$ defined by (6.8), where $V(m)$ is the irreducible 
representation of $U_q(sl_2)$ with highest weight $2m$. 

By definition, the trigonometric limit of $u_m(\l,\tau,\mu,p)$ 
is the leading coefficient
of the asymptotic expansion of $u_m(\l,\tau,\mu,p)$ as the modular
parameters $\tau,p$ tend to $i\infty$. We will denote this leading 
coefficient by $u_m(\l,\mu)$. 

Sending $\tau,p$ to $i\infty$ in the definition of \cite{FV2}, one obtains 
the following definition of the function $u_m(\l,\mu)$.

Let $q=e^t$ with $\text{Re}(t)>0$, and let $0<|A|<1$.
Define a function
$$
\gather
I_m(\l,\mu,A)= 
q^{-\l\mu}\int_{|T_1|=...=|T_m|=1}
\prod_{j=1}^m\frac{(T_jA^{-2}q^{\l+m}-q^{-\l-m})
(T_jA^{-2}q^{\mu+m}-q^{-\mu-m})}{(1-T_jA^2)(1-T_jA^{-2})}\times \\
\prod_{1\le i<j\le m}\frac{q^{-2}(1-T_iT_j^{-1})^2}
{(1-T_iT_j^{-1}q^{-2})(1-T_iT_j^{-1}q^2)}\wedge_{j=1}^m\frac{dT_j}
{2\pi \sqrt{-1}T_j}.\tag 8.1\endgather
$$ 
It is obvious that this function analytically 
continues to a rational function in $A$. We will denote this 
analytic continuation 
also by $I_m$. 

\proclaim{Definition} 
$$
u_m(\l,\mu)=I_m(\l,\mu,q^m).\tag 8.2
$$
\endproclaim

{\bf Remark.} Note that $A=q^m$ does not satisfy the condition $|A|<1$, 
which is why we needed to talk about the analytic continuation. 

The main result of this section is 

\proclaim{Theorem 8.1} Let $V(m)$ be the irreducible 
representation of $U_q(sl_2)$ with highest weight $2m$. Then  
$$
u_{V(m)}(\l,\mu)=q^{(3m-1)m}\frac{[2m]_q!}{m!}(q-q^{-1})^mu_m(\l,\mu),\tag 8.3
$$
where $u_{V(m)}$ is as in Section 6. 
\endproclaim

The rest of the section is the proof of Theorem 8.1. 

\proclaim{Lemma 8.2} 
$$
\int_{|T_1|=...=|T_m|=1}
\prod_{1\le i<j\le m}\frac{(1-T_iT_j^{-1})^2}
{(1-T_iT_j^{-1}q^{-2})(1-T_iT_j^{-1}q^2)}\wedge_{j=1}^m\frac{dT_j}
{2\pi \sqrt{-1}T_j}=q^{-m(m-1)/2}\frac{m!}{[m]_q!}.\tag 8.4
$$ 
\endproclaim

\demo{Proof} 
Let $0\le |p|<1$. 
Consider the Macdonald denominator of type $A_{m-1}$
(see \cite{M1}) 
$$
\Delta_{p,t}(x_1,...,x_m)=\prod_{i\ne j}\frac{(x_ix_j^{-1},p)_\infty}
{(tx_ix_j^{-1},p)_\infty},\tag 8.5
$$
where $(a,p)_\infty:=\prod_{j=0}^\infty(1-ap^j)$.

The Macdonald constant term identity (see \cite{M1}, p.20-21) 
says that the constant term 
of the Laurent series (8.5) (with respect to $x_i$) is 
given by
$$
c.t.(\Delta_{p,t})=m!\prod_{i<j}
\frac{(t^{j-i}p,p)_\infty(t^{j-i},p)_\infty}{(t^{j-i+1},p)_\infty
(t^{j-i-1}p,p)_\infty}.\tag 8.6
$$
Setting in this identity $p=0$, we get 
$$
c.t.(\Delta_{0,t})=m!\frac{(1-t)^m}{(1-t)...(1-t^m)}.\tag 8.7
$$
Substituting $t=q^{-2}$, we obtain the Lemma. 
$\square$\enddemo

Define the expression
$$
I_{k,m}=q^{-k(\l+\mu+2m)-\frac{k(k-1)}{2}}\frac{k!}{[k]_q!},\ 0\le k\le m,
\tag 8.8
$$ 
and the differential form
$$
\gather
\Omega_{k,m}=\prod_{j=1}^k
\frac{
(T_jA^{-2}q^{\l+m}-q^{-\l-m})
(T_jA^{-2}q^{\mu+m}-q^{-\mu-m})
(1-T_jA^{-2}q^{2m-2k-2})}
{(1-T_jA^2)
(1-T_jA^{-2}q^{2m-2k})
(1-T_jA^{-2}q^{-2})}\times\\ 
\prod_{1\le i<j\le k}
\frac{(1-T_iT_j^{-1})^2}
{(1-T_iT_j^{-1}q^{-2})(1-T_iT_j^{-1}q^2)}
\wedge_{j=1}^k\frac{dT_j}
{2\pi \sqrt{-1}T_j}.\tag 8.9\endgather
$$

\proclaim{Lemma 8.3} One has 
$$
\gather
\int_{|T_j|=A^2q^{-2(m-k)}(1+\epsilon)}\Omega_{k,m}
=\\
-kq^{-2(m-k)}\frac{(q^{\l+k}-q^{-\l-k})(q^{\mu+k}-q^{-\mu-k})(1-q^{-2})}
{(1-A^4q^{-2(m-k)})(1-q^{-2(m-k+1)})}
\int_{|T_j|=A^2q^{-2(m-k+1)}(1+\epsilon)}\Omega_{k-1,m}
+I_{k,m}.\tag 8.10
\endgather
$$
\endproclaim

\demo{Proof} Let us perform the integration with respect to $T_k$, for fixed 
$T_1,...,T_{k-1}$. 
It is obvious that 
the differential form $F_{k,m}$, as a function of $T_k$, 
 has two simple poles 
inside the circle of integration -- $T_k=A^2q^{2k-2m}$ and $T_k=0$. 
Therefore, the integral with respect to $T_k$ is equal to the 
sum of residues at these two poles. The residue at the first pole 
can be found by a direct computation and equals the first 
term on the right hand side of (8.10). 
The residue at zero equals to $I_{k,m}$ by Lemma 8.2. 
Lemma 8.3 is proved. 
$\square$\enddemo

Now let us prove Theorem 8.1. Let us move the contour of integration
in the definition of $I_m$ from $|T_j|=1$ to $|T_j|=A^2(1+\epsilon)$, 
via contours $|T_j|=B$, $A^2(1+\epsilon)\le B\le 1$.  
On the way, we do not run into any poles, therefore, we have
$$
u_m(\l,\mu)=q^{-\l\mu-m(m-1)}\int_{|T_j|=A^2(1+\epsilon)}\Omega_{m,m}|_{A=q^m}
\tag 8.11
$$
(in the sense of analytic continuation). 
So, to prove Theorem 8.1, it is enough to compute 
$\int_{|T_j|=A^2(1+\epsilon)}\Omega_{m,m}|_{A=q^m}$. 
We do it by using the recursive relation given in Lemma 8.3. 
Namely, we have 

\proclaim{Lemma 8.4} One has 
$$
\int_{|T_j|=A^2q^{-2(m-k)}(1+\epsilon)}\Omega_{k,m}|_{A=q^m}=
\sum_{j=0}^kc_{kj}I_{j,m},\tag 8.12
$$
where 
$$
c_{kj}=(-1)^{k-j}(1-q^{-2})^{k-j}\frac{k!}{j!}\prod_{i=j+1}^k
\frac{(q^{\l+i}-q^{-\l-i})(q^{\mu+i}-q^{-\mu-i})}
{(1-q^{2(i+m)})(1-q^{2(i-m-1)})}q^{2(i-m)}.\tag 8.13
$$
\endproclaim

\demo{Proof} The proof is a straightforward induction in $k$ 
using Lemma 8.3. 
$\square$\enddemo

Substituting $k=m$ in Lemma 8.4, and using the definition of 
$I_{k,m}$, we find the following expression for $u_m$:
$$
\gather
u_m(\l,\mu)=q^{-m(3m-1)}\frac{m!}{[2m]_q!}(q-q^{-1})^{-m}\times\\
q^{-\l\mu}\sum_{l=0}^m q^{-l(\l+\mu)-l(l-1)/2}(q-q^{-1})^l
\frac{[m+l]_q!}{[l]_q![m-l]_q!}\prod_{i=j+1}^m
(q^{\l+i}-q^{-\l-i})(q^{\mu+i}-q^{-\mu-i}).\tag 8.14
\endgather
$$
Comparing (8.14) with (7.12), we get Theorem 8.1. 

\head 9. Trace functions and Macdonald theory\endhead

In this section, 
following \cite{EK1} and 
\cite{FV1}, we will connect the results of this paper with 
the Macdonald-Ruijsenaars theory.

We restrict ourselves to the case of $\g=sl_n$,
$N=1$, and let $V$ be the q-analogue of the representation $S^{mn}\C^n$.
The zero-weight subspace of this representation is 1-dimensional, 
so the function $\Psi_V$ can be regarded as a scalar function. 
We will denote this scalar function by $\Psi_m(q,\l,\mu)$

Recall the definition of Macdonald operators \cite{M2,EK1}. 
They are operators on the space of functions 
$f(\l_1,...,\l_n)$ which are invariant under simultaneous shifting 
of the variables, $\l_i\to \l_i+c$, and have the form
$$
M_r=\sum_{I\subset \{1,...,n\}: |I|=r}
\left(\prod_{i\in I,j\notin I}
\frac{tq^{2\l_i}-t^{-1}q^{2\l_j}}{q^{2\l_i}-q^{2\l_j}}\right)T_I,\tag 9.1
$$  
where $T_I\l_j=\l_j$ if $j\notin I$ and $T_I\l_j=\l_j+1$ if 
$j\in I$. Here $q,t$ are parameters. We will assume that $t=q^{m+1}$, 
where $m$ is a nonnegative integer. 

It is known \cite{M2} that the operators $M_r$ commute. 
 From this it can be deduced that 
for a generic $\mu=(\mu_1,...,\mu_n)$, $\sum \mu_i=0$, there exists 
a unique power series $f_{m0}(q,\l,\mu)\in \C[[q^{\l_2-\l_1},...,
q^{\l_{n}-\l_{n-1}}]]$ such that the series
$f_m(q,\l,\mu):=q^{2(\l,\mu-m\rho)}f_{m0}(q,\l,\mu)$
satisfies difference equations 
$$
M_rf_m(q,\l,\mu)=
(\sum_{I\subset \{1,...,n\}: |I|=r} q^{2\sum_{i\in I}(\mu+\rho)_i})
f_m(q,\l,\mu).\tag 9.2
$$

{\bf Remark.}
The series $f_{m0}$ is convergent to an analytic (in fact, a trigonometric)
function.

The following theorem is contained in \cite{EK1}.

\proclaim{Theorem 9.1} (\cite{EK1}, Theorem 5) One has 
$$
f_m(q,\l,\mu)=\gamma_m(q,\l)^{-1}\Psi_m(q^{-1},-\l,\mu),\tag 9.3
$$
where 
$$
\gamma_m(q,\l):=
\prod_{i=1}^m\prod_{l<j}(q^{\l_l-\l_j}-q^{2i}q^{\l_j-\l_l}).\tag 9.4
$$
\endproclaim

{\bf Remark.} The exact statement of Theorem 5 of \cite{EK1}, in 
our conventions, is that the function 
$f_m(q,\l,\mu)\gamma_m(q,\l)$ is equal to 
$\text{Tr}|_{M_\mu}(\Phi_\mu^V(q^{-1})q^{2\l})$,
which is equivalent to Theorem 9.1.  

Let $\Cal D_W(q^{-1},-\l)$ denote the difference operator, obtained 
from the operator $\Cal D_W$ defined in Section 1 by the transformation 
$q\to q^{-1}$ and the change of coordinates $\l\to -\l$.
Let $\Lambda^r\C^n$ denote the q-analog of the r-th fundamental representation 
of $sl_n$. 
 
\proclaim{Corollary 9.2}
$$
\Cal D_{\Lambda^r\C^n}(q^{-1},-\l)=\delta_q(\l)\gamma_m(q,\l)\circ M_r
\circ \gamma_m(q,\l)^{-1}\delta_q(\l)^{-1}
$$
\endproclaim

\demo{Proof} This follows from Theorem 9.1 and Theorem 1.1. 
$\square$\enddemo

In conclusion 
of this section we would like to make several important remarks. 

{\bf Remark 1.} Corollary 9.2 is a degenerate (trigonometric) 
case of Theorem 5.2 in \cite{FV1}, 
which says that the elliptic Ruijsenaars operators are transfer 
matrices of the elliptic quantum $sl_n$ acting in $V[0]$. Thus,
Theorem 1.1 and Theorem 9.1
immediately imply the trigonometric case of Theorem 5.2 of \cite{FV1}
(i.e. the case without spectral parameter). 

{\bf Remark 2.} Conversely,  
the trigonometric case of Theorem 5.2 of \cite{FV1}
together with Theorem 1.1 immediately implies Theorem 9.1 (and many 
other results of \cite{EK1}). This is a ``direct'' proof 
of Theorem 9.1, in the sense that it involves
(unlike the original proof of \cite{EK1}) a direct 
computation of the radial parts of the central elements of $U_q(\g)$.
(Another direct proof of Theorem 9.1 is 
given in \cite{Mi}, where the radial part of the central element 
corresponding to the vector representation is computed).

{\bf Remark 3.} The line of 
argument discussed in Remark 2 can be extended to the elliptic case. Namely, 
combining an elliptic analogue of Theorem 1.1 (for affine Lie 
algebras at the critical level), and Theorem 5.2 of 
\cite{FV1}, one can prove an elliptic analogue of 
Theorem 9.1, which says that the radial parts of the central elements 
of $U_q(\widehat{sl_n})$ at the critical level corresponding to 
evaluation modules $\Lambda^r\C^n(z)$, 
acting on functions with values in $V[0]$, are elliptic Ruisjsenaars 
operators. This has been a conjecture for a number of years
(see e.g. \cite{Mi}, p.415). We plan to do this in a subsequent 
paper of this series.  

{\bf Remark 4.} In many arguments of this paper, Verma modules 
$M_\mu$ can be replaced with 
finite dimensional irreducible modules $L_\mu$ with sufficiently 
large highest weight, and one can prove analogs of Theorems
1.1-1.5 in this situation (in the same way). 
In particular, one may
set $\hat\Psi_m(q,\l,\mu)
=\text{Tr}(\hat\Phi_\mu^V q^{2\l})$, where \linebreak 
$\hat \Phi_\mu^V:L_\mu\to L_\mu\o V\o V^*[0]$ is the intertwiner with 
highest coefficient 1 (Such an operator exists iff $\mu-m\rho\ge 0$, 
see \cite{EK1}). Then one can show analogously to Theorem 9.1 
(see \cite{EK1}) that 
the function 
$\hat f_m(q,\l,\mu):=\gamma_m(q,\l)^{-1}\hat\Psi_m(q^{-1},-\l,\mu+m\rho)$ 
is the Macdonald polynomial $P_\mu(q,t,q^{2\l})$ with highest weight $\mu$
($\mu$ is a dominant integral weight). 
In this case, Theorem 1.1 says that Macdonald's polynomials 
are eigenfunctions of Macdonald's operators, Theorem 1.2 gives recursive 
relations for Macdonald's polynomials with respect to the weight 
(for $sl(2)$ -- the usual 3-term relation for orthogonal polynomials), 
and Theorem 1.3 is the Macdonald symmetry identity (see \cite{M2}).
(This representation theoretic derivation of the symmetry identity 
is somewhat different from the one in \cite{EK2}, where 
a pictorial argument is used.) 

\head 10. Limiting cases. \endhead

In this section we will discuss various degenerations of the function 
$F_{V_1,...,V_N}(q,\l,\mu)$, and the corresponding degenerate versions 
of Theorems 1.1-1.5. The main limiting cases we will be
interested in are the classical limit, and the rational limit.
The classical limit corresponds to passing
from $U_q(\g)$ to $\g$ in the trace construction; in this limit
the function $F$ depends rationally on $\mu$ but
trigonometrically of $\lambda$. This limit corresponds to 
the theory of spherical functions on the Lie group $G$ 
associated with $\g$, which is discussed 
in \cite{EFK1}. In the rational limit, which corresponds to the 
theory of spherical functions on $\g$ rather than $G$, 
the function $F$ becomes rational in both $\lambda$ and $\mu$,
restoring the symmetry. In this limit, the function $F$ is 
the Baker-Akhiezer function for a multivariable bispectral
problem (see \cite{Be}).

\subhead 10.1. The classical (KZB) limit\endsubhead

Let 
$$
F^{c}_{V_1,...,V_N}(\l,\mu)=\lim_{t\to 0}
F_{V_1,...,V_N}(q=e^t,\frac{\l}{2t},\mu).\tag 10.1
$$
We will call this limit the classical limit. 
The existence of this limit follows from 
Proposition 10.1 below. 

{\bf Remark.} Here and below we write the dependence of
functions on $q$ explicitly, since 
in this section $q$ is allowed to vary.  

{\bf Example 1.} If $\g=sl(2)$, N=1, and $V=V_1$ is the 
3-dimensional representation, we have 
$$
F^{c}_V(\l,\mu)=e^{-\l\mu/2}\frac{\mu}{\mu-1}
\left(1-\frac{1}{\mu}\frac{1+e^\l}{1-e^\l}\right).\tag 10.2
$$

The classical limit is obtained when in the situation 
of Section 1 we take the ordinary enveloping algebra $U(\g)$ 
instead of the quantized one $U_q(\g)$. 

More precisely, let $\Phi_\mu^V$ be intertwining operators for $U(\g)$ defined 
as in Section 1, and set 
$$
\Psi^c_{V_1,...,V_N}(\l,\mu)=
\text{Tr}((\Phi^{V_1}_{\mu+\sum_{i=2}^Nh^{(*i)}}\o 1^{N-1})
...\Phi^{V_N}_\mu e^{\l}),\tag 10.3
$$
Also, set 
$$
\delta(\l)=e^{(\l,\rho)}\prod_{\alpha>0}(1-e^{-(\l,\alpha)}),\tag 10.4
$$
and let $Q^c(\mu)$ be the limit of $\Q(\mu)$ as $q\to 1$ 
(i.e. it is defined as in Section 1 from representation theory of $U(\g)$).
Then we have 

\proclaim{Proposition 10.1} 
$$
F^c_{V_1,...,V_N}(\l,\mu)=\delta(\l)
[\Q^{-1}(\mu)^{(*N)}\o ...\o \Q^{-1}(\mu-h^{(*2...*N)})^{(*1)}]
\Psi^c_{V_1,...,V_N}(\l,-\mu-\rho).\tag 10.5
$$
\endproclaim

\demo{Proof} The proof is straightforward.
\enddemo

Let us now look at the degenerations of the properties of $F_{V_1,...,V_N}$ 
in the classical limit. We start with the analogue of Theorem 1.1. 

First of all, we have the following analogue of Proposition 2.1, 
which is proved analogously to Proposition 2.1.  
 
\proclaim{Proposition 10.2} (i) For any  
element $X$ of $U(\g)$ there exists a unique 
differential operator $D_X$ acting on $V[0]$-valued functions, 
such that 
$$
\Tr(\Phi^V_\mu Xe^{\l})=D_X\Tr(\Phi^V_\mu e^{\l}).\tag 10.6
$$

(ii) If $X$ is central then $D_{XY}=D_YD_X$
for all $Y\in U(\g)$. In particular, if $X$, 
$Y$ are central then $D_XD_Y=D_YD_X$.  
\endproclaim

The operator $D_X$ can be computed explicitly for any element $X$, but 
in general the answer is complicated. However, if $X$ is the quadratic 
Casimir $C$, the answer is easy to write down.
Namely, define $\tilde D_X=\delta(\l) D_X\delta(\l)^{-1}$. 
Then we have (\cite{E},\cite{ES1}):
$$
\tilde D_C=\Delta_\h-\sum_{\alpha>0}
\frac{f_\alpha e_\alpha}{2\text{sinh}^2\frac{1}{2}(\l,\alpha)}
-(\rho,\rho),\tag 10.7
$$
where $f_\alpha,e_\alpha$ are root generators such that 
$(e_\alpha,f_\alpha)=1$, and $\Delta_\h$ is the Laplacian on the 
Cartan subalgebra associated with the standard invariant form.

Thus, we have the following classical analogue of 
Theorem 1.1:

\proclaim{Theorem 10.3} For any $X$ in the center of $U(\g)$, let 
$p_X$ be the symmetric polynomial on $\h^*$ such that 
$X|_{M_\mu}=p_X(\mu+\rho)$. Then we have 
$$
\tilde D_X^\l F^c_{V_1,...,V_N}(\l,\mu)=p_X(-\mu)F^c_{V_1,...,V_N}(\l,\mu).
\tag 10.8
$$
In particular, 
$$
(\Delta_\h-\sum_{\alpha>0}
\frac{f_\alpha e_\alpha}{2\ \text{sinh}^2\frac{1}{2}(\l,\alpha)})
F^c_{V_1,...,V_N}(\l,\mu)=(\mu,\mu)F^c_{V_1...,V_N}(\l,\mu).\tag 10.9
$$
\endproclaim

Formula (10.9) was obtained in \cite{E,ES1}, but it can also be  
 derived by taking the classical limit in Theorem 1.1. 

Now let us consider the classical analogue of Theorem 1.2. 
Let $\Cal D_W^{\vee,c}$ denote the 
difference operators defined by formula (1.12) 
for $q=1$ 
(i.e. $\RR(\mu)$ are the exchange matrices for $U(\g)$ with 
$\mu$ replaced by $-\mu-\rho$). Then we have 
the following result, obtained by passing to the limit in 
Theorem 1.2. 

\proclaim{Theorem 10.4} 
$$
\Cal D_W^{\vee,c,\mu} F^c_{V_1...V_N}(\l,\mu)=\chi_W(e^\l)
F^c_{V_1...V_N}(\l,\mu),
\tag 10.10
$$
\endproclaim

{\bf Example 2.} 
In the case of Example 1,
Theorems 10.3 and 10.4 have the form
$$
\left(\frac{\d^2}{\d \l^2}-\frac{1}{2\ \text{sinh}^2(\l/2)}\right)
F^c_V(\l,\mu)=\frac{\mu^2}{4}F^c_V(\l,\mu),
$$
and
$$
(T+\frac{(\mu-2)(\mu+1)}{\mu(\mu-1)}T^{-1})F^c_V(\l,\mu)=
(e^{\l/2}+e^{-\l/2})F^c_V(\l,\mu) 
$$
(where $T$ is the shift by 1 in $\mu$), which is easily checked 
from (10.2).  

Now consider the classical limit of Theorem 1.3. 
For this purpose introduce the classical dynamical 
r-matrix $r(\l)$, which is the classical limit of 
the exchange matrix $R(q,\l)$. This matrix is defined by the formula
$$
R(q=e^t,\frac{\l}{2t})=1-2r(\l)t+O(t^2),\tag 10.11
$$
and equals to 
$$
r(\l)=-\frac{1}{2}\Omega+\frac{1}{2}
\sum_{\alpha>0}\text{cotanh}\frac{1}{2}(\l,\alpha)
e_\alpha\wedge f_\alpha,
$$
where $\Omega$ is the Casimir tensor (see \cite{EV}). 
Taking the quasiclassical limit in 
Theorem 1.3, and using that $r(-\l)=r^{21}(\l)$,
we obtain the following result. 

\proclaim{Theorem 10.5} For any $j=1,...,N$, one has
$$  
\gather
\biggl[\frac{\d}{\d h^{(j)}}-
(\sum_{l<j} r_{lj}(\l)-\sum_{l>j} r_{jl}(\l))\biggr]F^c_{V_1,...,V_N}
(\l,\mu)=\\
[(\mu+\frac{1}{2}\sum x_i^2)_{*j}
+\sum_{l<j}\sum x_i^{(*j)}\o x_i^{(*l)}]F^c_{V_1,...,V_N}(\l,\mu),
\tag 10.12
\endgather
$$
where $\frac{\d}{\d h^{(j)}}X(\l)=\frac{\d}{\d \nu}X(\l)$ 
if $X$ is a tensor-valued function whose $j$-th component has weight $\nu$.  
\endproclaim

The last equation is the trigonometric limit of the KZB equation, which is why 
the classical limit is called ``the KZB limit''.

Let us now consider the classical limit of Theorem 1.4.
Let $K_j^{\vee,c}$ be the difference operators defined 
by formula (1.17) for $q=1$ 
(i.e. $\RR(\mu)$ are exchange matrices for $U(\g)$ with 
$\mu$ replaced by $-\mu-\rho$). Then we have 
the following result, obtained by passing to the limit in 
Theorem 1.4.

\proclaim{Theorem 10.6} One has
$$
K_j^{\vee,c}F^c_{V_1...V_N}(\l,\mu)=
(e^\l)_jF^c_{V_1...V_N}(\l,\mu).\tag 10.13
$$
\endproclaim

Finally, Theorem 1.5 does not have an analogue in the classical limit. 
In this limit, the symmetry between $\l$ and $\mu$ is destroyed, since
$F^c$ is a product of $e^{(\l,\mu)}$ 
with a function that is trigonometric in $\l$ but rational in $\mu$. 
 
\subhead 10.2. The rational limit\endsubhead

The rational limit is a further degeneration of the classical limit. Namely, 
let 
$$
F^{r}_{V_1,...,V_N}(\l,\mu)=\lim_{\gamma\to 0}
F^c_{V_1,...,V_N}(\l\gamma,\mu/\gamma).\tag 10.14
$$
We will call this limit the rational limit. 
The existence of this limit and the fact that 
$det(F^r)\ne 0$ can be deduced from Corollary 3.3 
in \cite{ES1}.

{\bf Example 3.} If $\g=sl(2)$, N=1, and $V=V_1$ is the 
3-dimensional representation, we have 
$$
F^{r}_V(\l,\mu)=e^{-\l\mu/2}
\left(1+\frac{2}{\l\mu}\right).\tag 10.15
$$

The degeneration of Theorem 1.1 
in this limit is the following theorem. 
Let $\tilde D_X^r$ be the rational limit of 
$\tilde D_X$, i.e. $\tilde D_X^r(\l)$
is the leading coefficinet of $\tilde D_X(\gamma\l)$
as $\gamma\to 0$. For instance, 
$$
\tilde D_C^r=
\Delta_\h-\sum_{\alpha>0}
\frac{2f_\alpha e_\alpha}{(\l,\alpha)^2}.\tag 10.16
$$

\proclaim{Theorem 10.7} For any $X$ in the center of $U(\g)$, let 
$p_X$ be the symmetric polynomial on $\h^*$ such that 
$X|_{M_\mu}=p_X(\mu+\rho)$. 
Let $p_X^r$ be the top degree component of $p_X$. Then we have 
$$
\tilde D_X^{r,\l} 
F^r_{V_1,...,V_N}(\l,\mu)=p_X^r(-\mu)F^r_{V_1,...,V_N}(\l,\mu).
\tag 10.17
$$
In particular, 
$$
(\Delta_\h-\sum_{\alpha>0}
\frac{2f_\alpha e_\alpha}{(\l,\alpha)^2})
F^r_{V_1,...,V_N}(\l,\mu)=(\mu,\mu)F^r_{V_1...,V_N}(\l,\mu).\tag 10.18
$$
\endproclaim

The degeneration of Theorem 1.2 looks as follows:

\proclaim{Theorem 10.8} Equations (10.17),(10.18) are satisfied for 
the function $F^{r,*}_{V_N^*,...,V_1^*}$. 
\endproclaim

{\bf Example 4.} In the situation of Example 3, 
Theorems 10.7, 10.8 have the form 
$$
\left(\frac{\d^2}{\d \l^2}-\frac{2}{\l^2}\right)F^r_V(\l,\mu)=
\frac{\mu^2}{4}F_V^r(\l,\mu), 
$$
$$
\left(\frac{\d^2}{\d \mu^2}-\frac{2}{\mu^2}\right)F^r_V(\l,\mu)=
\frac{\l^2}{4}F_V^r(\l,\mu), 
$$
which is easily checked from (10.15). 

Using the asymptotics of $F^r_{V_1,...,V_N}(\l,\mu)$ at infinity, 
similarly to arguments of Section 5, one can deduce 
from Theorems 10.7,10.8 the following analog of 
Theorem 1.5 (the symmetry theorem):

\proclaim{Theorem 10.9}  
The function $F^r_{V_1...V_N}$ is symmetric: 
$$
F^r_{V_1...V_N}(\l,\mu)=F^{r,*}_{V_N^*...V_1^*}(\mu,\l),\tag 10.19
$$
\endproclaim

Thus, the symmetry, lost in the first limit, is restored after taking 
the second limit. 

{\bf Remark 1.} Another proof of Theorem 10.8 is based on representation 
of the above sequence of two limits as a single limiting procedure, 
which is symmetric in $\l$ and $\mu$. Namely, one can show that 
$$
F^r_{V_1,...,V_N}(\l,\mu)=\lim_{s,t\to 0}
F_{V_1,...,V_N}(q=e^{st/2},\frac{\l}{t},\frac{\mu}{s}),\tag 10.20
$$
after which Theorem 10.9 follows from Theorem 1.5. 

{\bf Remark 2.} Theorems 10.7, 10.8 show that 
the function $F^r_{V_1,...,V_N}(\l,\mu)$ is a solution of the matrix 
{\it bispectral 
problem} in several variables (on the bispectral 
problem, see e.g. \cite{DG,G}). The Baker-Akhiezer function of the rational 
Calogero system of type $A$, 
which is a known solution of the 
multidimensional bispectral problem (\cite{VSC}, see also \cite{Be}), 
is a special case of $F^r_{V_1,...,V_N}(\l,\mu)$ ($N=1$, 
$V=V_1=S^{mn}\C^n$).      

Finally, let us consider the rational limit of Theorems 1.3, 1.4. 
To formulate the analog of Theorem 1.3, introduce the 
rational limit of the classical dynamical r-matrix, 
$r^0(\l)=\lim_{\gamma\to 0}\gamma r(\gamma\l)$. It has the form 
$$
r^0(\l)=
\sum_{\alpha>0}
\frac{e_\alpha\wedge f_\alpha}{(\l,\alpha)}.\tag 10.21
$$
Taking the rational limit in 
Theorem 10.4, we get 

\proclaim{Theorem 10.10}
For any $j=1,...,N$, one has
$$  
\gather
\biggl[\frac{\d}{\d h^{(j)}}-
(\sum_{l<j} r^0_{lj}(\l)-\sum_{l>j} r^0_{jl}(\l))\biggr]F^r_{V_1,...,V_N}
(\l,\mu)=\\
\mu_{*j}F^r_{V_1,...,V_N}(\l,\mu).
\tag 10.22
\endgather
$$
\endproclaim

The analogue of Theorem 1.4 is 

\proclaim{Theorem 10.11} Equation (10.22) is satisfied for 
the function $F^{r,*}_{V_N^*,...,V_1^*}$. 
\endproclaim

\subhead 10.3. The qKZ and KZ limits\endsubhead

Assume that $|q|<1$. 
The qKZ limit is defined by 
$$
F^{qKZ}_{V_1,...,V_N}(\l,\mu)=\lim_{(\l,\alpha_i)\to -\infty}
q^{2(\l,\mu)}F_{V_1,...,V_N}(\l,\mu).\tag 10.23
$$
It is easy to check using Theorem 50 of \cite{EV} that 
$$
\gather
F^{qKZ}_{V_1,...,V_N}(\l,\mu)=\\
[\Q^{-1}(\mu)^{(*N)}\o ...\o \Q^{-1}(\mu-h^{(*2...*N)})^{(*1)}]
\<(\Phi^{V_1}_{\mu+\sum_{i=2}^Nh^{(*i)}}\o 1^{N-1})
...\Phi^{V_N}_\mu \>=\\
[\Q^{-1}(\mu)^{(*N)}\o ...\o \Q^{-1}(\mu-h^{(*2...*N)})^{(*1)}]
\JJ^{1...N}(\mu)^*,\tag 10.24\endgather
$$
where $\<,\>$ denotes the highest matrix element. 
(The last expression is an endomorphism 
of $(V_N^*\o...\o V_1^*)[0]$, which is regarded as an element of 
\linebreak $(V_1\o...\o V_N)[0]\o (V_N^*\o...\o V_1^*)[0]$). 
In particular, this function is independent on $\l$. 

Let us now consider the behavior of the 
equations given by Theorems 1.1-1.5 in the 
qKZ limit.

 The MR equations given by Theorem 1.1 become trivial. 
Namely, when $(\l,\alpha_i)\to -\infty$, one has $\JJ(\l)\to 1$, 
and hence $\RR(\l)\to \Cal R^{21}$. The matrix $\Cal R^{21}$ is triangular, 
so only its diagonal part contributes to the trace. Inspection of this
diagonal part shows that $\lim_{(\l,\alpha_i)\to -\infty}
\Cal D_W=\sum_\nu \text{dim}W[\nu]\ T_\nu$, and 
the limiting equation is 
$$
(\sum_\nu \text{dim}W[\nu]\ q^{-2(\nu,\mu)}T_\nu)F^{qKZ}=
\chi_W(q^{-2\mu})F^{qKZ}
$$
(where $T_\nu$ is the shift of $\l$),
which is a trivial consequence of the fact that 
$F^{qKZ}$ is independent on $\l$. 

 The dual MR equations given by Theorem 1.2 have a 
slightly more interesting limit.
It is easy to see that the only term on each side of (1.11) which survives
in the limit is the term corresponding to the lowest weight 
$\nu_W$ of $W$. Therefore, the limiting equation has the form
$$
(v_W^*\o 1,\RR _{WV_N^*}^{01}(\mu+h^{(*1...*N-1)})...
\RR_{WV_1^*}^{0N}(\mu)(v_W\o 1))
F_{V_1,...,V_N}^{qKZ}(\mu+\nu_W)= 
F_{V_1,...,V_N}^{qKZ}(\mu).\tag 10.25
$$
 
The qKZB equations become the trigonometric limit of the qKZ equations.
Namely, for all $j=1,...,N$ we have 
$$
[\Cal R_{j,j+1}^{-1}...\Cal R_{jN}^{-1}
(q^{-2\mu})_j\Cal R _{1j}...\Cal R_{j-1,j}
\o D_j]F^{qKZ}_{V_1,..,V_N}(\mu)=
F^{qKZ}_{V_1,..,V_N}(\mu).\tag 10.26
$$
If $N=2$, these equations are closely related to 
the Arnaudon-Buffenoir-Ragoucy-Roche (ABRR) equation (see Lemma 2.4). 
In general, they are essentially the N-component version of the 
ABRR equation. 

{\bf Remark.} The dual qKZB equations do not seem to have 
a reasonable qKZ limit. Also, 
the symmetry relation (Theorem 1.5) does not hold in the qKZ limit since 
the function depends on $\mu$ and not on $\lambda$. 

The KZ limit is obtained from the qKZ limit as $q\to 1$, in which case the 
qKZ equations degenerate into the trigonometric KZ equations 
(the quasiclassical limit of the ABRR equation). We leave it to the reader
to derive the limiting equations in this case. 

We plan to consider these limits in more detail in another paper in 
the more interesting case of affine Lie algebras and quantum affine algebras.

\end